%% file: ricchiuto_etal.tex
\documentclass[11pt]{article} 

\usepackage[T1]{fontenc} 							
\usepackage[utf8]{inputenc} 

\usepackage{type1cm}        
\usepackage{makeidx}         
\usepackage{graphicx}        
\usepackage{multicol}        
\usepackage[bottom]{footmisc}
\usepackage{newtxtext}       
\usepackage{newtxmath}       

\usepackage{epsfig}
\usepackage{hyphenat} 
\usepackage{subfig}
\usepackage{tikz}
\usepackage{pgfplotstable}
\usepackage{filecontents,comment}
\usepackage{algorithm}
\usepackage{algcompatible}
\usepackage[noend]{algpseudocode}
\newcommand{\mario}[1]{\textcolor{black}{#1}}

\makeindex  

\definecolor{ForestGreen}{HTML}{228b22}
\definecolor{Emerald}{HTML}{50C878}
\definecolor{NavyBlue}{HTML}{000080}
\definecolor{Orange}{HTML}{ffa500}
\definecolor{BrickRed}{HTML}{CB4154}
\include{newcommands}

\pgfplotstableset{
    create on use/H/.style={create col/copy column from table={h.txt}{0}}}
\pgfplotstableset{
    create on use/Hsquare/.style={create col/copy column from table={hsquare.txt}{0}}}
\pgfplotstableset{
    create on use/Htau/.style={create col/copy column from table={htau.txt}{0}}}
\usetikzlibrary{spy}
           
\begin{document}

\title{$h-$ and $r-$ adaptation on simplicial meshes using  MMG tools}
\author{Luca Arpaia$^{(1)}$, H\'elo{\"\i}se Beaugendre$^{(2)}$, Luca Cirrottola$^{(2)}$, Algiane Froehly$^{(2)}$, \\Marco Lorini$^{(2)}$, L\'eo Nouveau$^{(3)}$ and Mario Ricchiuto$^{(2)}$\\[15pt] 
%
%
%
$^{(1)}$  Coastal Risk and Climate Change Unit, French Geological Survey, \\3 Av. C. Guillermin 45060 Orl\'eans Cedex 2,  France\\[5pt]  
$^{(2)}$   INRIA, Universit\'e de Bordeaux, CNRS, Bordeaux INP, IMB UMR 5251,  \\200 Avenue de la Vieille Tour, 33405 Talence cedex, France, \\[5pt]
$^{(3)}$  INSA Rennes, CNRS, IRMAR - UMR 6625, F-35000 Rennes, France}


\maketitle

%

\abstract{We review some recent work on the enhancement and application of both $r-$ and $h-$  adaptation  techniques,
benefitting of the functionalities of  the remeshing  platform Mmg:  {\tt www.mmgtools.org}.   Several contributions revolve around 
the level-set adaptation capabilities of the platform. These have been used   to identify  complex surfaces and  then to either produce
conformal 3D meshes, or to define a metric allowing to perform $h$-adaptation and control geometrical errors in the context of  immersed boundary  flow simulations.
The performance of the recent distributed memory parallel   implementation ParMmg is also discussed.
In a similar spirit, we  propose some  improvements of $r-$adaptation methods to handle embedded fronts.} 


\section{Introduction}
\label{sec:general}
\input{hrintro.tex}

\section{$h-$adaptation: embedded geometries, parallel implementation} 
\label{sec:hgeneral}
\input{hgeneral.tex}

\subsection{Level-sets and IBM with $h-$adaptation}
\label{sec:ibmresults}
\input{hresults.tex}

\subsection{Parallel $h-$adaptation}
\label{sec:hparallel}
\input{parmmg.tex}

\section{$r-$adaptation for    embedded geometries and fronts} 
\label{sec:rgeneral}
\input{rgeneral.tex}

\input{rresults2D.tex} 

\section{Summary and perspectives}
\label{sec:conclusions}
This paper reviews recent  work performed using on the capabilities of the  remeshing software platform  \texttt{Mmg}.  
As other  open source packages for computational mechanics or mesh generation, this platform  allows 
to apply its existing functionalities to problems of interest, as well as developing new meshing/remeshing methods.

Ongoing work related to the  results presented are related to the robustification mesh deformation method in order to account for
general curved boundaries  \cite{CIRROTTOLA2021110177},
to the investigation of the impact of adaptation in the framework of higher order  embedded methods~\cite{Scovazzi3,Scovazzi4},
to the full support of  non-manifold surface adaptation and anisotropic remeshing in parallel, to the exploration of hybrid parallelization strategies.

%
%

\section*{Acknowledgements}
Dedicated to our dear friend C\'ecile Dobrzynski.

\bibliographystyle{spmpsci}
\bibliography{biblio1,biblio,LibMesh}

\end{document}

%% file: newcommands.tex
\renewcommand{\(}{\ensuremath{\begin{pmatrix}}}  
\renewcommand{\)}{\ensuremath{\end{pmatrix}}}

\newcommand{\bbm}{\ensuremath{\begin{bmatrix}}}
\newcommand{\ebm}{\ensuremath{\end{bmatrix}}}
\newcommand{\bpm}{\ensuremath{\begin{pmatrix}}}
\newcommand{\epm}{\ensuremath{\end{pmatrix}}}

\newcommand{\rst}[1]{\ensuremath{{\mathbin\upharpoonright}\raise-.5ex\hbox{$#1$}}}

\newcommand{\be}{\begin{equation}}
\newcommand{\ee}{\end{equation}}

\newcommand{\ba}{\begin{aligned}}
\newcommand{\ea}{\end{aligned}}

\newcommand{\grad}{\boldsymbol{\nabla}}

\newcommand{\nbf}{\mathbf{n}}

\newcommand{\xbf}{\mathbf{x}}

\newcommand{\wbf}{\mathbf{w}}

\newcommand{\zerobf}{\mathbf{0}}
\newcommand{\deltabf}{\boldsymbol{\delta}}

\newcommand{\chibf}{\boldsymbol{\chi}}

\newcommand{\epsilonbf}{\boldsymbol{\epsilon}}

\newcommand{\xibf}{\boldsymbol{\xi}}

\newcommand{\sigmabf}{\boldsymbol{\sigma}}

\newcommand{\Qcal}{\ensuremath{\mathcal{Q}}}

\newcommand{\SketchScale}{1.0\textwidth}

%% file: hrintro.tex
Mesh adaptation has become nowadays   a powerful tool to improve the discrete representation of complex solution fields in many applications, and particularly in  computational fluid dynamics~\cite{Park2016}. 
Adapting the mesh  may  lead to a non-negligible computational overhead, as well as to more complex algorithmic and software developments. 
This motivates the quest for efficient and robust methods.

$h$-adaptation allows to optimize the discrete representation of a field of interest by inserting and removing mesh entities. 
This has proven to be a very powerful tool, and it is today  quite    mature and generic (cf. \cite{Park2016} and references therein). 
Several  aspects of this approach are still quite complex in general. This is the case for conservative high order 
solution transfer between  meshes with different topologies, which is a non-trivial step with a non-negligible  cost ~\cite{FarrellMaddison2011,Alauzet2016,hermes2018highorder,KucharikShashkov2008,Re2017}.
$h$-adaptation also entails  non foreseeable, non-linear, and dynamic  changes in the distribution of the computational workload,  which is an inherently weakly parallelizable
feature on distributed memory architectures.

Conversely,   node relocation without topology change (or $r-$adaptation) has been the first adaptation approach for structured grids, later extended to unstructured ones.
It remains  still alluring due to the possibility of a minimally intrusive coupling with existing computational solvers, with no modification of  the  data structures. 
Moreover, devising conservative projections is quite natural with $r-$adaptation, due to the inherently continuous nature of the process. This allows to exploit
space-time or Arbitrary-Lagrangian-Eulerian techniques to build high order  conservative remapping~\cite{Kucharik2003efficient,Trulio1961,ThomasLombard1979}
compliant  with the Geometric Conservation Law (GCL).\\

The results of this paper benefit of the functionalities of the Mmg platform
\cite{MmgPlat}     which implements some well known   research on metric-based mesh adaptation~\cite{Dapogny2014Mmg,Dobrzynski2008Anisotropic} in  a free and open source software application and library for simplicial mesh adaptation~\cite{MmgSW}.  The platform  provides adaptive and implicit domain remeshing through three libraries and applications \texttt{Mmg2d}, \texttt{Mmgs} and \texttt{Mmg3d}, targeting two-dimensional cartesian and surfacic triangulations, and three-dimensional tetrahedral meshes. Recent efforts are devoted to the development of its parallel version    \texttt{ParMmg}~\cite{CirrottolaFroehly2019,ParMmgSW} on top of the \texttt{Mmg3d} remesher.
A moving mesh library \texttt{Fmg}, which uses the \texttt{Mmg} library for its mesh data structure and geometry representation, is under development. 
Example of applications in various fields and documentation   can be found on the platform website~\cite{MmgPlat},
also providing tutorials  under the section \textit{Try Mmg!}, and library examples  in the respective repositories~\cite{MmgSW,ParMmgSW}. 

The contribution here discussed  exploit the  functionalities of the platform.  We present the work
done in the context of the application of  $h$-adaptation to  control geometrical errors in   immersed boundary  flow simulations, and we discuss the performance of the parallel   implementation in ParMmg.
In a similar spirit,  we review  some  improvements and applications of $r-$adaptation methods to flows involving immersed solid  boundaries, and moving flronts.
For immersed boundaries we discuss some ideas used to track embedded fronts and level sets.  \\

%
%
%
%
%

 \mario{The paper is organized in two main parts. Sect.~\ref{sec:hgeneral} is devoted to the application and implementation  of $h-$adaptation methods, with focus on 
 level-set adaptation, application to immersed boundary methods, and parallel $h-$adaptation.
 Sect.~\ref{sec:rgeneral} presents instead advances and applications of $r-$adaptation methods with focus on the resolution of embedded fronts. 
The paper is ended by an overlook on ongoing activities.}

%% file: hgeneral.tex
Anisotropic $h-$adaptation has largely proved useful in the context of numerical simulations to capture the physical behavior of complex phenomena at a reasonable computational cost~\cite{FreyAlauzet2005}.
Following a now classical approach, we consider adaptation strategies based on a local error estimation, which is converted in a map for  the size and orientation of the mesh edges.
This is done via an iterative process based on successive  flow field evaluations, error estimations, and a-posteriori local mesh modifications.
As described in~\cite{FreyGeorge2008}, mesh sizes and edge directions are controlled via 
 metric tensors   built starting from error estimations involving   the Hessian of the  target output. 
The eigenvalues $\lambda_i$ of these matrices are directly linked to the sizes $h_i$ of the elements edges in the directions $i$ (where $\lambda_i = 1/h_i^2$), with these directions given by the eigenvectors.
Different criteria can be used to evaluate metrics, two examples relevant for the applications discussed later are  recalled hereafter.

If several metric fields are available, e.g. the ones of Sects.~\ref{sec:hphysical} and \ref{sec:hlevelset}, one may seek to combine them.
To this end, we use here the simultaneous reduction method described in~\cite{FreyGeorge2008}.
Denoting by $\mathcal{M}_1$ and $\mathcal{M}_2$ the two metrics defined at the same node, the resulting metric $\mathcal{M}_{1 \cap 2}$ respects both sizes prescribed by the two initial metrics.
With a geometrical analogy, the ellipsoid $\mathcal{E}_{1 \cap 2}$ associated to $\mathcal{M}_{1 \cap 2}$ is the biggest ellipsoid included in both $\mathcal{E}_{1}$ and $\mathcal{E}_{2}$, associated to $\mathcal{M}_1$ and $\mathcal{M}_2$, respectively.

The metrics for physical and level-set adaptation will be recalled in Sects.~\ref{sec:hphysical} and \ref{sec:hlevelset}. Then, an application case of level-set adaptation to immersed boundary simulations will be shown in Sect.~\ref{sec:ibmresults} and an example of level-set discretization will be presented in Sect. \ref{sec:sagrada}
Finally, we will shortly discuss in Sect.~\ref{sec:hparallel} our extension of h-adaptation techniques to parallel computing environments.

\subsection{Physical adaptation}
\label{sec:hphysical}

If we aim at controlling the error between a certain output field, $u$, and its linear interpolation, we can exploit
the upper bound on a mesh element $K$~\cite{FreyAlauzet2005}:

\begin{equation}
|| u - \Pi_h u ||_{\infty,K} \leq C_d \max_{\mathbf{e} \in K} \langle \mathbf{e}, \mathcal{M}_1(K) \mathbf{e} \rangle\;.
\end{equation}

Here, $\mathbf{e}$ denotes the edges of the mesh element and $C_d$ is a constant depending on the dimension.  The metric $\mathcal{M}_1(K)$ is a function of  the Hessian of $u$, $\mathcal{H}(u)$:

\begin{equation}
\mathcal{M}_1 = {^t}\mathcal{R} \begin{pmatrix}
    \lambda_1       & 0 & 0 & \\
    0       & \lambda_2 & 0 & \\
    0       & 0 & \lambda_3 & 
\end{pmatrix} \mathcal{R}\;,
\end{equation}

\noindent with $\mathcal{R}$ the matrix of the eigenvectors of $\mathcal{H}(u)$ and with $\lambda_i$ defined as:

\begin{equation}
\lambda_i = \min \left( \max \left( |h_i|, \frac{1}{h_{max}^2} \right), \frac{1}{h_{min}^2} \right)\;,
\end{equation}

\noindent where $h_i$ is the $i$-th eigenvalue of the Hessian,   $h_{min}$ (resp. $h_{max}$) are  the minimum (resp. maximum) allowed sizes for the mesh edges. 

\subsection{Level-set adaptation}
\label{sec:hlevelset}

The principle is here to improve the accuracy in the definition of the zero iso-value of a level-set  function,  $\Phi$, describing the distance from a given surface.
This allows to adapt  the volume mesh on which the function is defined.  To this end, we can use the metric field proposed in ~\cite{Dapogny2012},  and given by:

\begin{equation}
\mathcal{M}_2 = {^t}\mathcal{R} \begin{pmatrix}
    \frac{1}{\epsilon^2}       & 0 & 0 & \\
    0       & \frac{\left| \lambda_1 \right|}{\epsilon} & 0 & \\
    0       & 0 & \frac{\left| \lambda_2 \right|}{\epsilon} & 
\end{pmatrix} \mathcal{R}\;,
\end{equation}

\noindent where $\mathcal{R}=(\nabla \Phi \;v_1 \;v_2)$, with $(v_1,\;v_2)$  a basis of the tangent plane to the local iso-value surface of $\Phi$, $\lambda_i$  the eigenvalues of its Hessian and $\epsilon$ a user-defined target approximation  error. 
As proposed in the above reference, this  metric field can be imposed  in a  user-defined  layer of thickness  $w$,  close to   the zero iso-level.
Outside this region,   the mesh size is set to grow  linearly up to $h_{max}$, unless  other constraints  are set.

%% file: hresults.tex
Unfitted discretizations are becoming quite popular for the additional geometrical flexibility they offer.
In these methods,  the accuracy with which   boundary conditions are met is  directly linked to the accuracy of the implicit description of the solid~\cite{MittalIaccarino2005}.
Mesh adaptation with respect to the level-set  is a useful tool  to control this error. We exploit this in the context of a   Brinkman-type volumic penalization discretization of the Navier-Stokes equations~\cite{KhadraAngot2000}.
This method falls in the category of continuous forcing immersed boundary methods, meaning that the flow equations are solved throughout the whole domain and the enforcement of rigid motion inside the immersed object is done via a continuous volumic forcing term.
In the case of the Navier-Stokes equations, this leads to a compact form

\begin{equation}
  \frac{\partial \mathbf{u}}{\partial t}
  +\nabla \cdot \mathbf{F}_C(\mathbf{u})
  +\nabla \cdot \mathbf{F}_V(\mathbf{u},\nabla\mathbf{u})
  +\mathbf{p}(\chi, \mathbf{u}) = {\bf 0}\;\;,\label{NScompact}
\end{equation}

\noindent where $\mathbf{u} \in \mathbb{R}^{m}$ denotes the vector of the $m$ conservative variables, $\mathbf{p}\in\mathbb{R}^{m}$ the penalization term, $d$ the space dimension, $\mathbf{F}_C, \mathbf{F}_V\in\mathbb{R}^{m}\otimes\mathbb{R}^d$ the inviscid and viscous flux functions.
The localization of the immersed object within the computational Navier-Stokes domain is done through a mask function, $\chi$, which is the Heaviside function of the level-set representation of the immersed solid.

\begin{figure}[htbp]
\hspace{-0.5cm}
\begin{minipage}{0.4\textwidth}
\includegraphics[bb= 0in 0in 11.69in 8.27in, width = \linewidth]{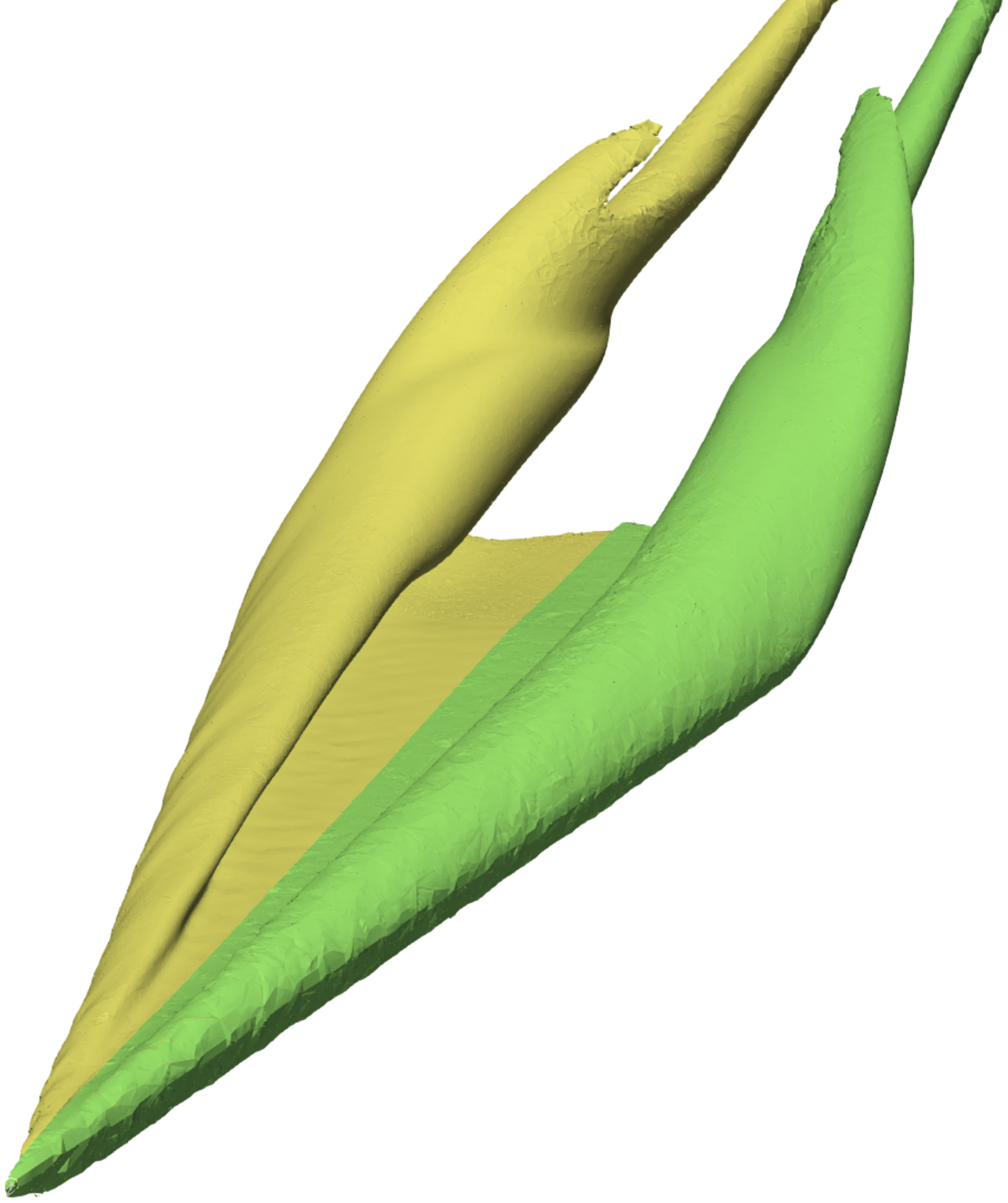}
\end{minipage}
\begin{minipage}{0.65\textwidth}
\includegraphics[width = 0.495\linewidth,trim=400 250 400 150,clip]{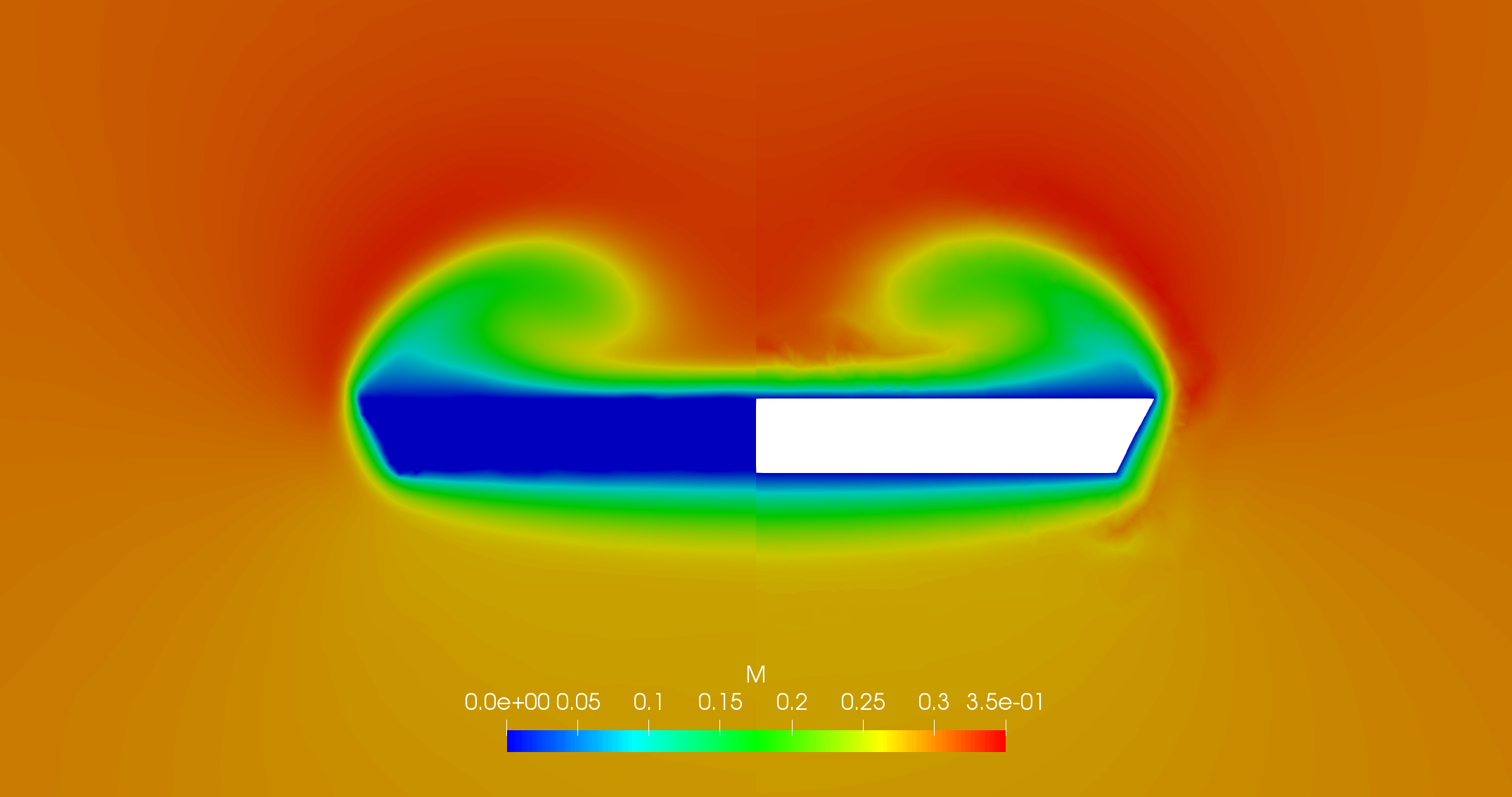} \includegraphics[width = 0.495\linewidth,trim=400 300 400 100,clip]{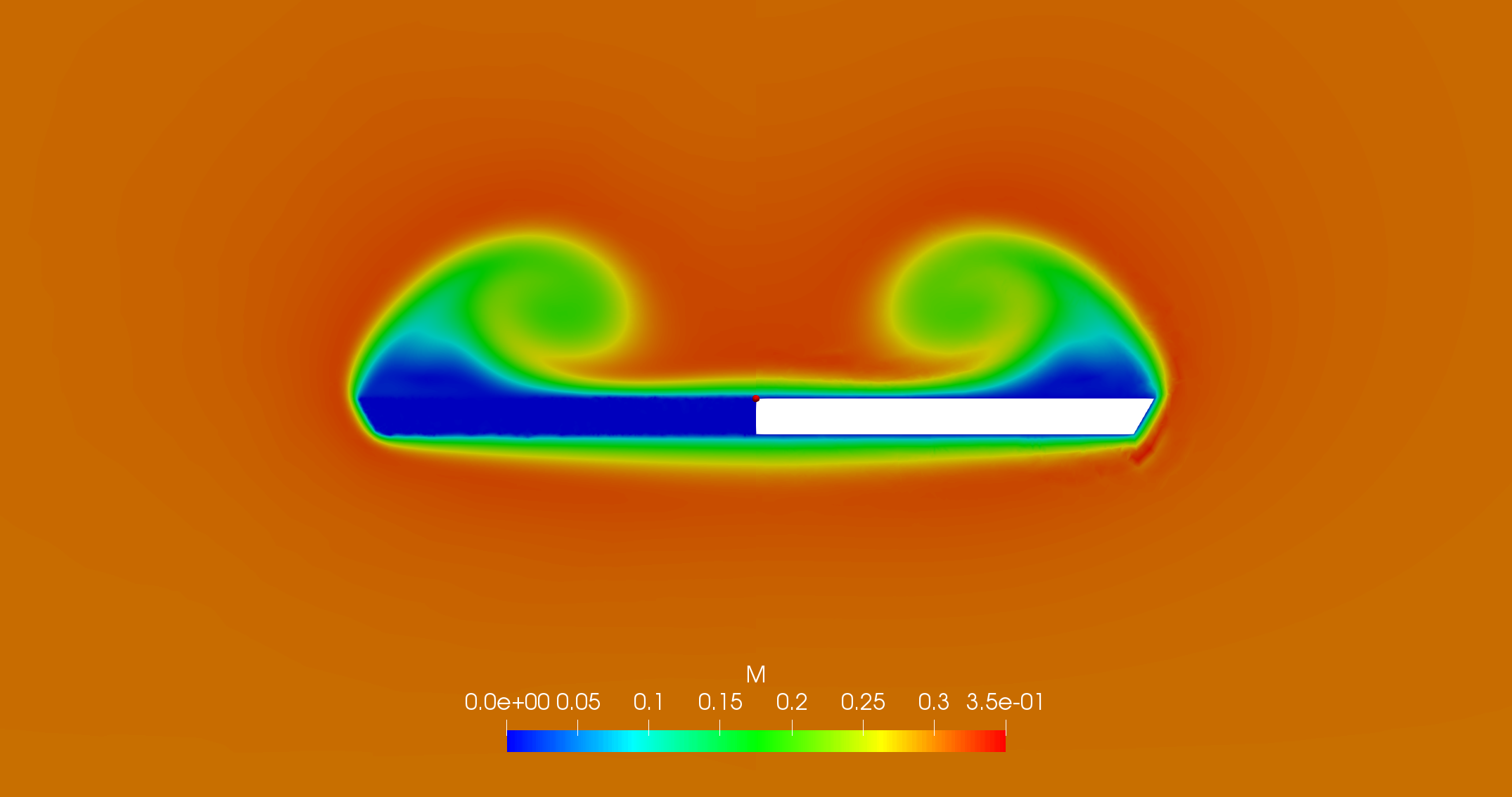}\\
 \includegraphics[width = 0.495\linewidth,trim=400 400 400 100,clip]{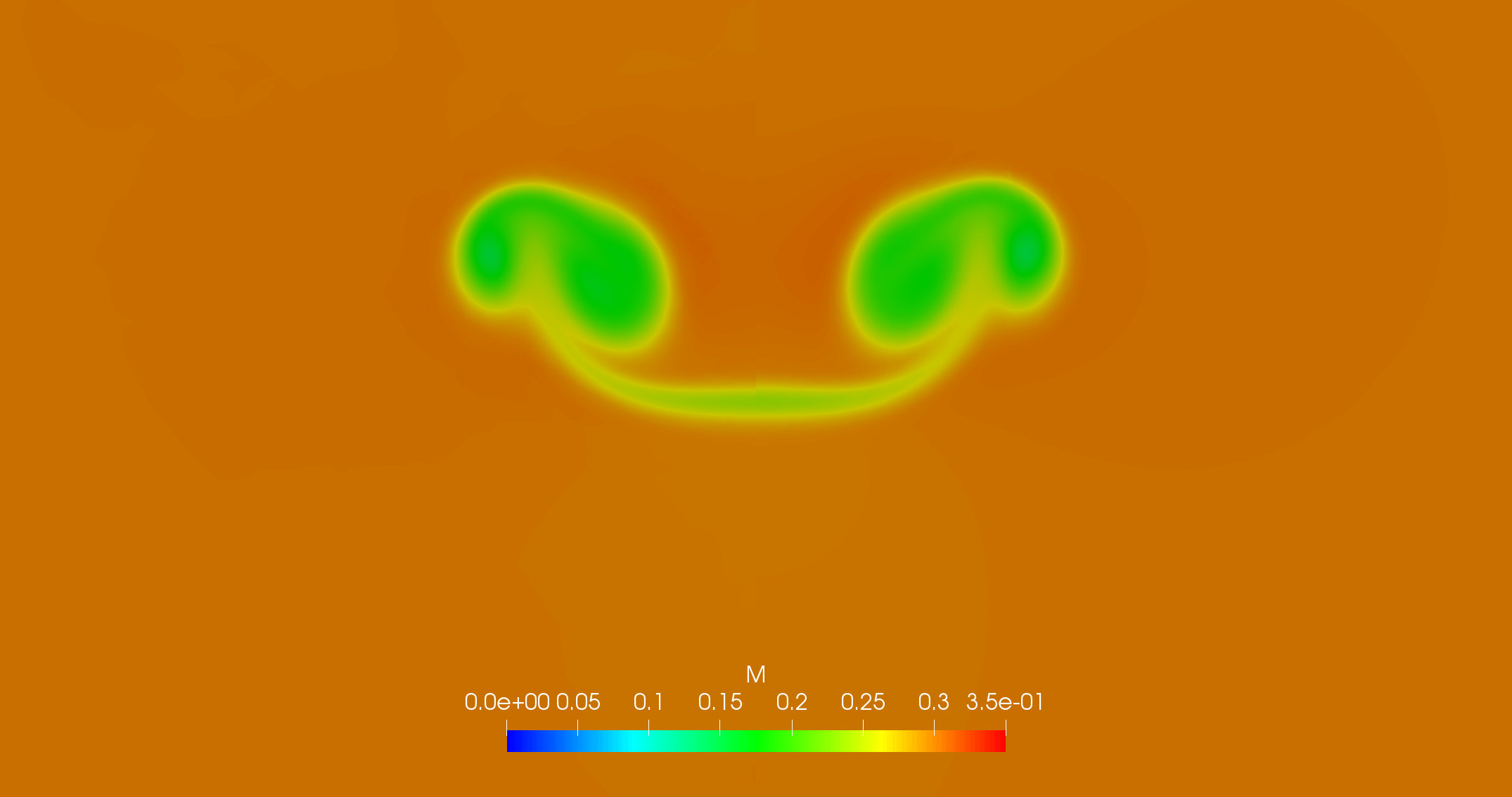} \includegraphics[width = 0.495\linewidth,trim=400 450 400 50,clip]{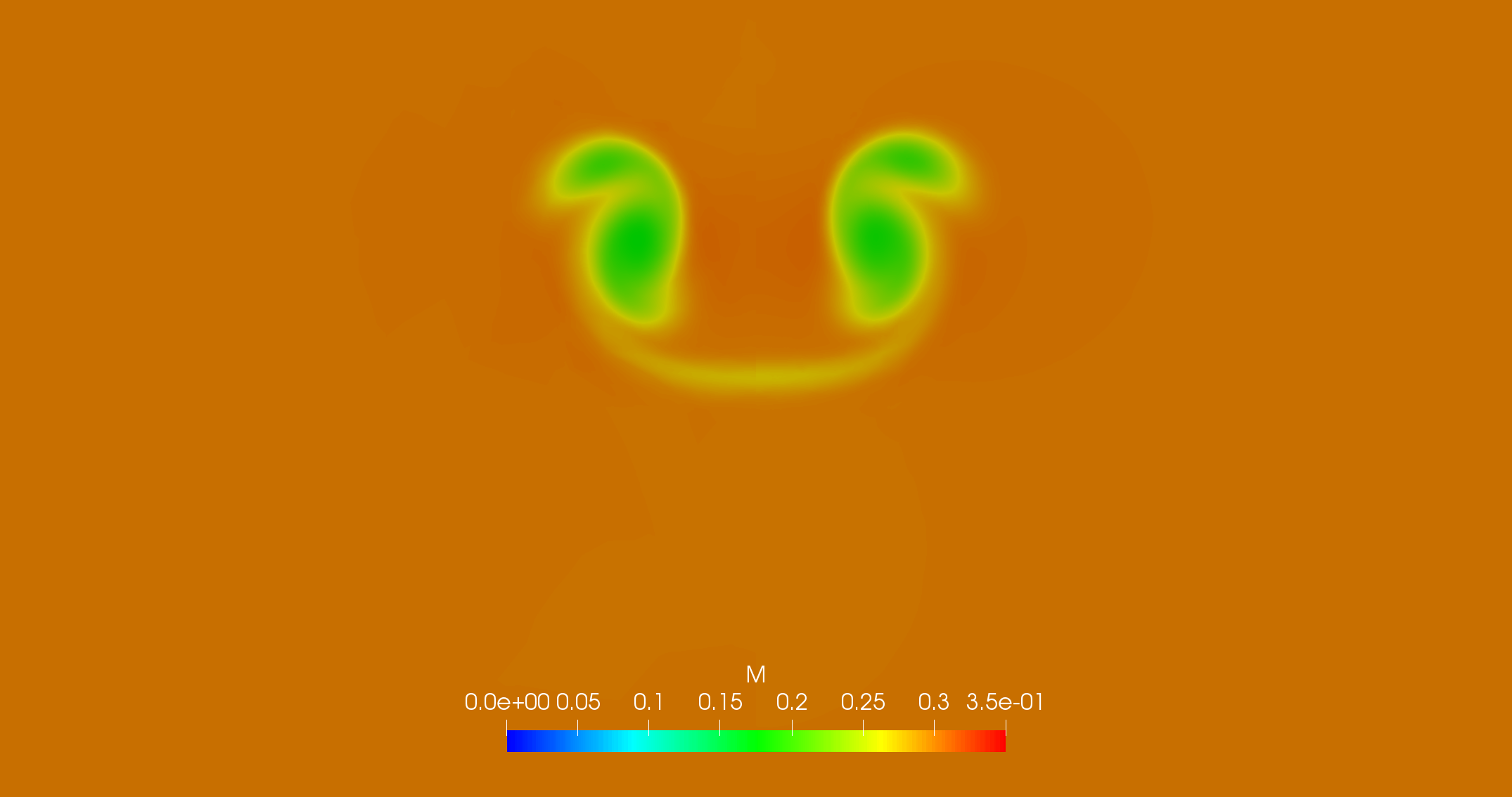}
\end{minipage}\caption{Laminar delta wing. \textbf{left} Roll-up of primary and secondary vortices. \textbf{right} $\mathbb{P}^2$ Mach number distribution at $x=0.5C$, $x=C$, $x=1.5C$ and $x=2C$. In each picture on the  left half 
 the  unfitted computation,  on the  right the  body-fitted one.}\label{fig:wing_adaptation} 
\end{figure}

We present in the following results obtained by coupling $h-$adaptation with a \mario{nodal discontinnous Galerkin discretization of the immersed boundary model of Eq.~\ref{NScompact}.  The 
discretization choices are discussed  in~\cite{Lorini2018Eccomas}.}  The test case considered is a steady laminar flow at high angle of attack around a delta wing with a sharp leading edge. 
This is a benchmark test case for adaptive methods for vortex dominated external flows.
Even if the geometry is not complex as the one in Sect.~\ref{sec:sagrada}, $h-$adaptation can be important for the imposition of the immersed boundary condition in order to correctly reproduce the flow topology.
See for an example Fig.~\ref{fig:wing_adaptation}, where the primary and secondary vortices and the DG $\mathbb{P}^2$ Mach number distribution at $x=0.5C$, $x=C$, $x=1.5C$ and $x=2C$ ($C$ being the chord of the wing) are compared between the immersed boundary (left half) and a body-fitted simulation (right half).
The mesh for the unfitted simulation consists of 446294 tetrahedra and is adapted to both the Mach number distribution and the level-set (cf. Sects.~\ref{sec:hphysical} and \ref{sec:hlevelset}).

\begin{figure}[htbp]
    \centering    
    \includegraphics[width=0.49\textwidth,trim=0mm 0mm 0mm 0mm,clip]{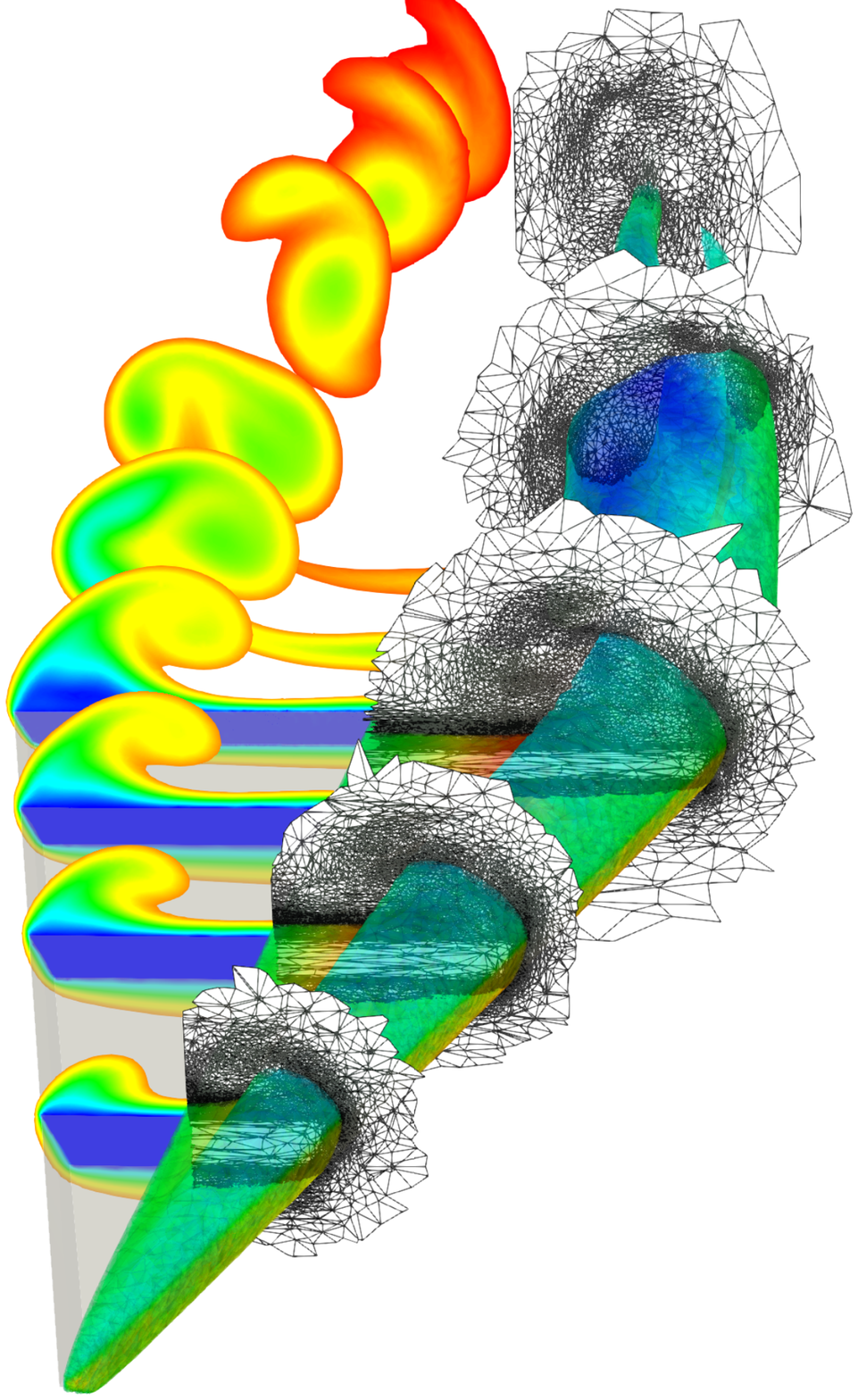}   
    \raisebox{-0.0125\height}{\includegraphics[width=0.49\textwidth,trim=0mm 0mm 0mm 0mm,clip]{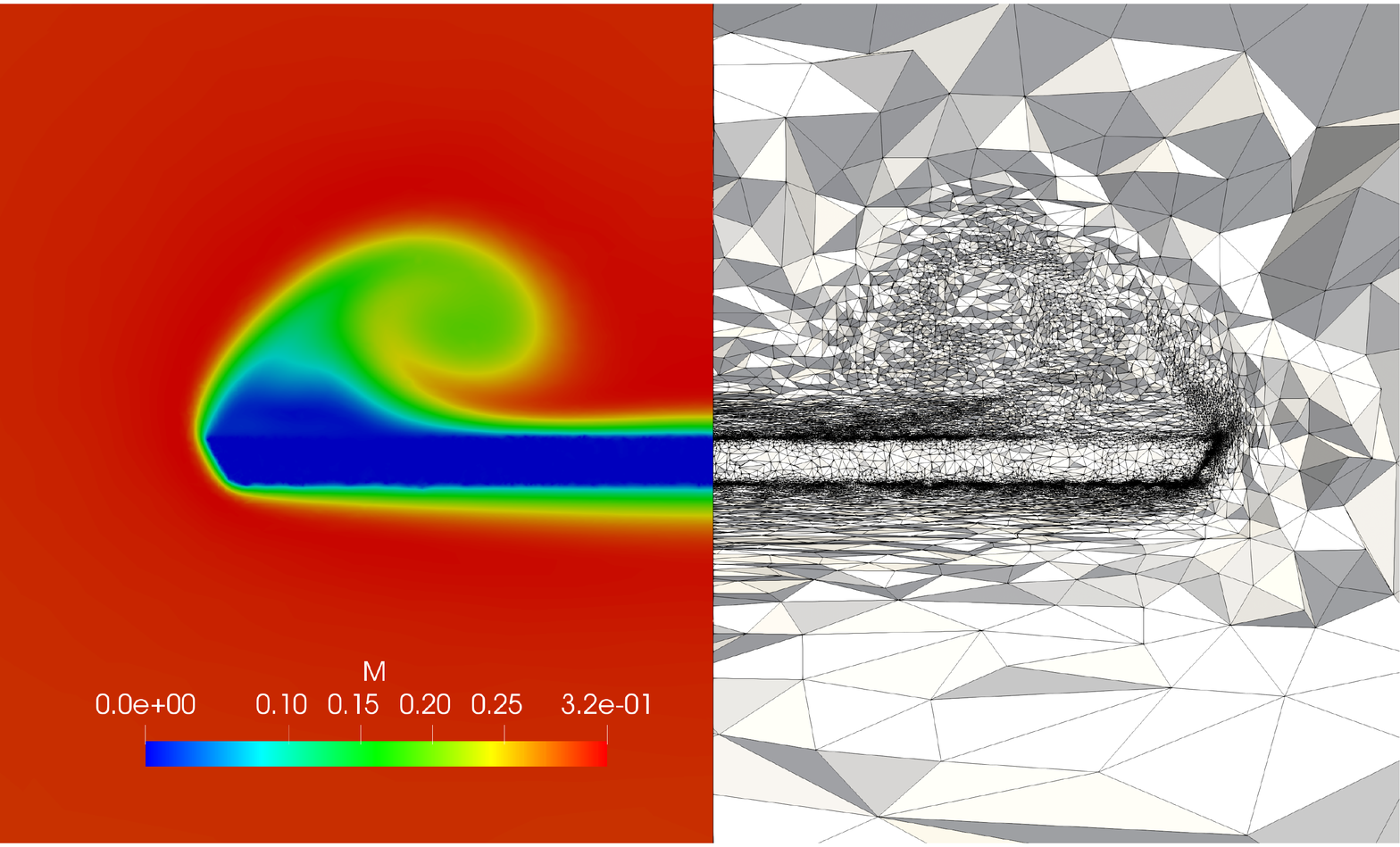}}
    \caption{Laminar delta wing. \textbf{left} Slices of the Mach number and adapted mesh. \textbf{right} Comparison flow - mesh at the trailing edge.}\label{fig:deltawing}  
\end{figure}

Then, as a proof of concept, starting from the DG $\mathbb{P}^2$ simulation on the adapted mesh, the full wing mesh has been adapted to fine level of detail (28 milions tetrahedra) with respect to a metric that takes into account both the level-set and the computed Mach number distribution.
In Fig.~\ref{fig:deltawing}(left) the left half of the wing presents some slices of the Mach number distribution ($\mathbb{P}^2$ approximation) at different wing locations and in the wake. 
The surface of the wing is represented only as a reference for the reader.
The right half of the wing presents a Mach number isosurface showing the correct representation of the flow topology (with the primary and secondary vortices) and some slices of the adapted mesh.
Fig.~\ref{fig:deltawing}(right) shows a slice of the adapted mesh at the trailing edge of the wing compared with the DG solution. 
The differences between the imposed metrics, that closely follow the Mach distribution and the level-set description of the wing, can be clearly seen in the different regions of the domain.
Finally, in Table~\ref{table:deltawing} we present the mesh statistics in terms of quality of the edges with respect to the imposed metrics from the original to the final mesh. 

\begin{table}[!t]
\begin{center}
\caption{Laminar delta wing, mesh statistics for $h-$adaptation from Original to Final mesh. Percentage of edges $ N^{(a,b]} $ whose length in the assigned metrics falls in the interval  $ l_{\mathcal{M}} \in (a,b] $.}
\label{table:deltawing}
\begin{tabular}{p{1.7cm}p{1cm}p{1.1cm}p{1.2cm}p{1.2cm}p{1.2cm}p{1.2cm}p{1.1cm}p{1cm}}
\hline\noalign{\smallskip}
Mesh (\# tets) & $ N^{(0   , 0.3 ]}  $ &
    $ N^{(0.3 , 0.6 ]} $ &
    $ N^{(0.6 , 0.7 ]} $ &
    $ N^{(0.71, 0.9 ]} $ &
    $ N^{(0.9 , 1.3 ]} $ &
    $ N^{(1.3 , 1.41]} $ &
    $ N^{(1.41, 2   ]} $ &
    $ N^{> 2         } $ \\
\noalign{\smallskip}\hline\noalign{\smallskip}
O (446294) & 0.93\% & 2.81\% & 1.20\% & 2.01\% & 2.99\% & 0.64\% & 3.19\% & 86.23\%\\ 
F (28083948) & 0.01\% & 0.58\% & 3.20\% & 24.36\% & 63.18\% & 6.22\% & 2.45\% & 0\%\\
\noalign{\smallskip}\hline\noalign{\smallskip}
\end{tabular}
\end{center}
\end{table}

\subsection{An example of implicit domain meshing: La Sagrada Familia}\label{sec:sagrada}
Level-set methods represent the boundary of a closed body as a level-set (typically the zero level-set) of a continuous function. This point of view can be particularly useful when trying to generate the volume mesh of a body starting only from an approximate representation of its boundary, for example a surface triangulation from an STL file, not necessarily displaying the sufficient regularity necessary for volume mesh generation (for example, being intersecting or non-orientable). In this case, this starting triangulation can be embedded in a larger volume mesh and a signed distance function from the surface triangulation can be defined on mesh nodes.
The mesh is  adapted according to the metrics defined in the previous section. The procedure can  be iterated to  increase the geometrical accuracy of the model.
A new surface triangulation can also be generated in correspondence of the zero level set, allowing to extract the interior body volume mesh.

\begin{figure}[htbp]
    \centering
    \subfloat[Level set function.]{
      \includegraphics[width=0.35\textwidth,trim=0mm 0mm 0mm 0mm,clip]{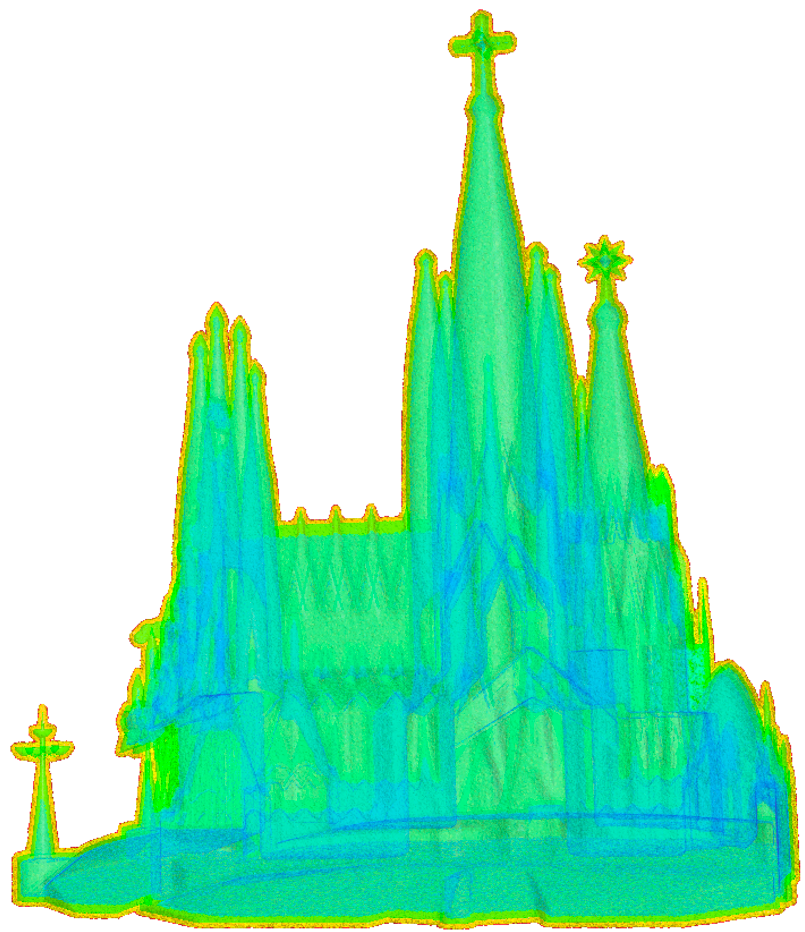}\label{fig:sagradaLS}}
    \subfloat[Volume remeshing.]{
      \includegraphics[width=0.55\textwidth,trim=0mm 0mm 0mm 0mm,clip]{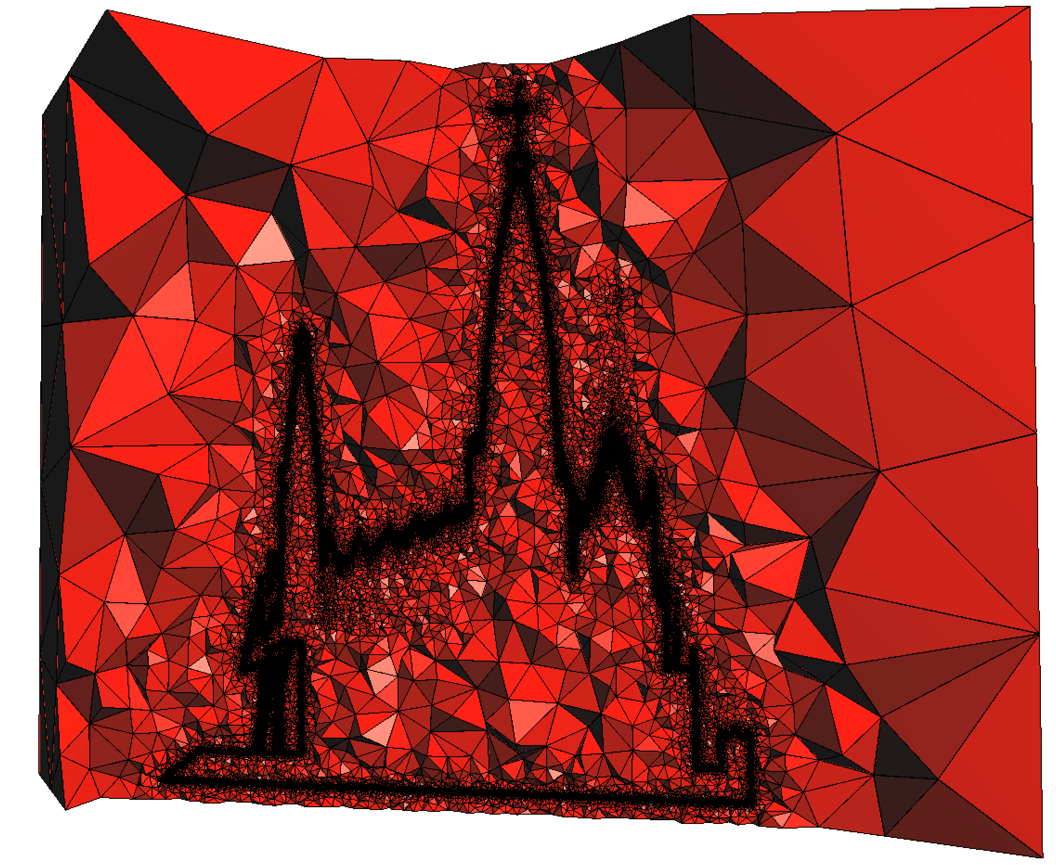}\label{fig:sagradaVol}}\\
    \subfloat[Surface discretization.]{
      \includegraphics[width=0.4\textwidth,trim=0mm 0mm 0mm 0mm,clip]{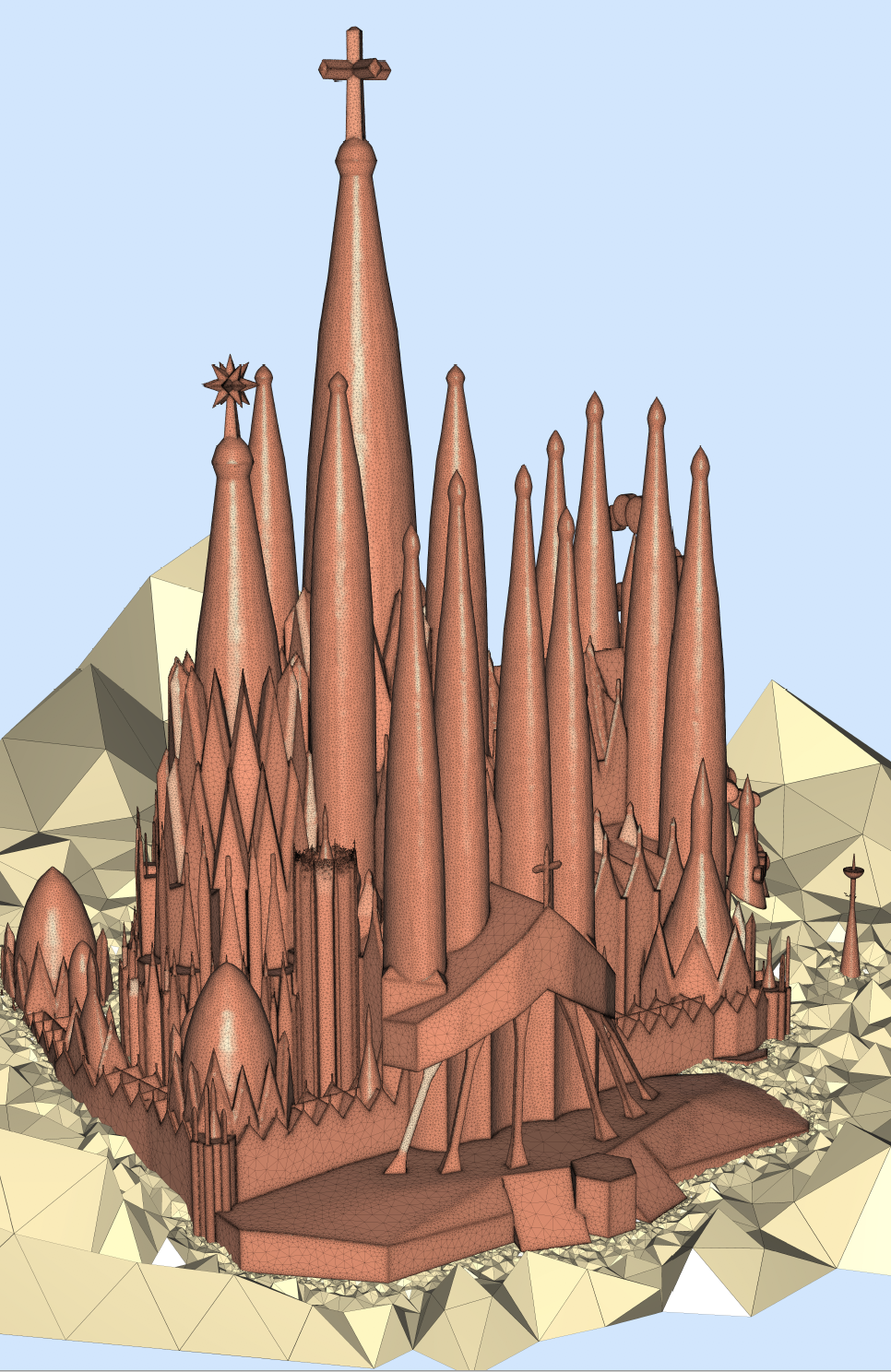}\label{fig:sagradaSurf}}
    \subfloat[Detail of the volume remeshing.]{
      \includegraphics[width=0.5\textwidth,trim=0mm 0mm 0mm 0mm,clip]{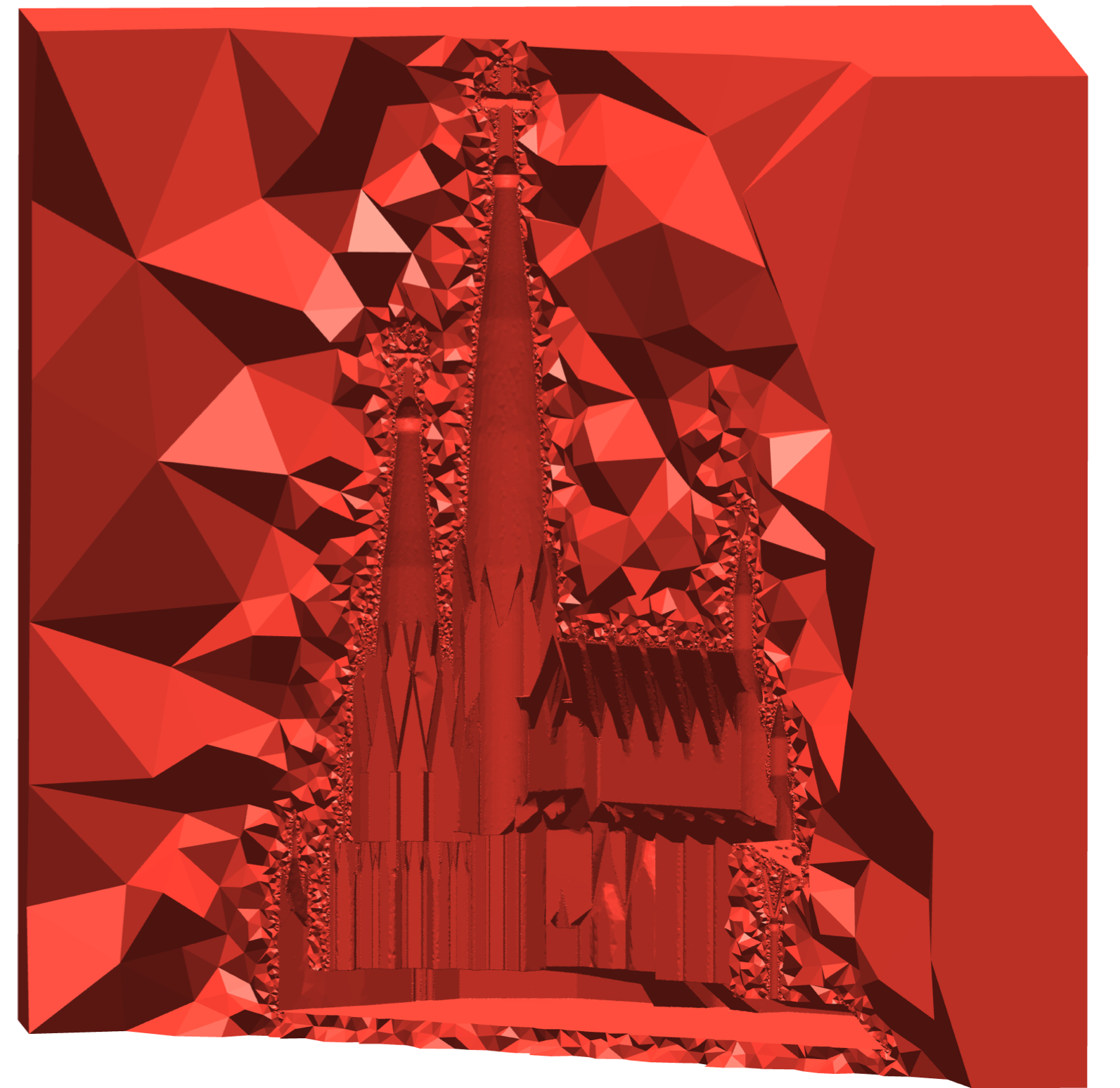}\label{fig:sagradaVolCut}}\\
    \caption{Implicit domain remeshing and surface discretization of the Sagrada Familia building with \texttt{Mmg3d}. From 50 to 150 million tetrahedra on a desktop computer, with 4 remeshing iterations.}\label{fig:sagrada}
\end{figure}

An application example is given in Fig.~\ref{fig:sagrada}, where the implicit domain meshing procedure is applied to the Sagrada Familia starting from the STL files provided by the International Meshing Roundtable for the meshing contest of the 2017 edition\footnote{https://imr.sandia.gov/26imr/MeshingContest.html},
using \texttt{Mmg}~\cite{Dapogny2014Mmg,MmgSW} for the implicit domain remeshing and \texttt{Mshdist}~\cite{Dapogny2012,MshdistSW} for the signed distance function computation.
{\color{black}
Isotropic remeshing is selected for this case}, and four remeshing iterations are performed. The initial mesh has 9 million nodes and 51 million tetrahedra, and the final one has 26 million nodes and 153 million tetrahedra (the worst tetrahedron has isotropic quality\footnote{The isotropic length of an edge AB is defined as \[ l_{AB} = \frac{||AB||}{h_A - h_B}\log{\frac{h_B}{h_A}} \] while the isotropic tetrahedron quality is defined as \[ Q = \alpha \frac{V}{(\sum_{i=1}^{6}l_j^2)^{1/3}} \] with  $ \alpha $ a  normalization factor to get quality equal to 1 for the regular tetrahedron with unit edges.}
equal to 0.2, while 99\% of the tetrahedra have quality above 0.5, and 95\% of the edges have unit-lengths between 0.7 and 1.4).

%% file: parmmg.tex
Computational mechanics solvers nowadays routinely exploit parallel, distributed memory computer architectures, raising the need for generating and adapting meshes whose size, in terms of computer memory, is larger and larger.
Even if the distributed mesh is to be handled by a sequential remesher by gathering it on a single process, it is possible that a single computing node cannot store it in memory due to its memory limitation. Also, sequential remeshing in a parallel simulation represents a significant performance bottleneck~\cite{Park2016AIAA}. Parallel remeshing is thus becoming increasingly demanded in large scale simulations.

The \texttt{ParMmg} application and library~\cite{ParMmgSW} is built on top of the \texttt{Mmg3d} remesher to provide parallel tetrahedral mesh adaptation in a free and open source software. Among the many possible remeshing parallelization methods (for example~\cite{Castanos1997}~\cite{DeCougny1999}~\cite{Oliker2000}~\cite{Chrisochoides2003}~\cite{Flaherty1998Parallel}~\cite{Cavallo2005Parallel}~\cite{Digonnet2017Massively}~\cite{Bernard2016MeshAdapt}), a modular approach is adopted by selecting an iterative remeshing-repartitioning scheme that does not modify the adaptation kernel~\cite{Bernard2016MeshAdapt}. As described in depth in~\cite{CirrottolaFroehly2019}, the sequential remeshing kernel is applied at each iteration on the interior partition of each process while maintaining fixed (non-adapted) parallel interfaces. Then the adapted mesh is repartitioned in order to move the non-adapted frontiers to the interior of the partitions at the next iteration, thus progressively eliminating the presence of non-adapted zones as the iterations progress. The repartitioning algorithm currently explicitly displace the parallel interfaces by a given number of element layers, before migrating the mesh parts generated by the displacement.

A weak scaling test is performed by refining a sphere of radius 10 while keeping the workload of each process as constant as possible as the number of processes is increased. The test is performed on the \texttt{bora} nodes of the \texttt{PlaFRIM} cluster\footnote{{\tt www.plafrim.fr}}, and the input and output mesh data are shown in Table~\ref{tab:weakTable}. The weak scaling performances are shown  on the left in Fig.~\ref{fig:scaling}. The slow, steady increase in the time spent in the repartitioning and redistribution phase shows that there is still space for optimizing the parallel mesh redistribution scheme.

A strong scaling test is then performed with an isotropic sizemap $ h $ describing a double Archimedean spiral
\begin{equation}
  h(x,y) = \min( 1.6 + | \rho - a \theta_1 | + 0.005, 1.6 + | \rho + a \theta_2 | + 0.0125 )
\end{equation}
with
\begin{equation}
  \begin{tabular}{ll}
      $ \theta_1 = \phi + \pi \left( 1+\text{floor} \left( \frac{\rho}{2 \pi a} \right) \right) $ \\
      $ \theta_2 = \phi - \pi \left( 1+\text{floor} \left( \frac{\rho}{2 \pi a} \right) \right) $
    \end{tabular}
\end{equation}
and $ \phi = \text{atan2}(y,x) $, $ \rho = s \sqrt{ x^2 + y^2 }$, $ a = \text{0.6} $, $ s = \text{0.5} $, into a sphere of radius 10 with uniform unit edge length. The surfacic adaptation and a volumic cut of the adapted meshes on 1 and 1024 processes are shown in Fig.~\ref{fig:Archimede}.
The resulting edge length statistics are presented in Table~\ref{tab:strongTable}, {\color{black}
showing the variation in the distribution of the edge lengths as the number of processes is increased.
It can be noticed that there is a slight trend for an increase of large edge sizes with more processes. This can be explained by the fact that, when increasing the number of partitions on the same initial mesh, the number of interfaces also increases, requiring more work to the remesher to refine the volume mesh near the coarse parallel interfaces, which are frozen by the algorithm, and making the argument for tuning the number of remeshing iterations (six in this test) with the number of processes. This trend can be also visually appreciated in Fig.~\ref{fig:scaling}.
The time performances in log scale are shown  on the right in Fig.~\ref{fig:scaling}.
} The speedup $ S_p $ over $ p $ processes is defined as the ratio between the wall time on 1 process and the wall time on $ p $ processes
\begin{equation}
  S_p =  T_1/T_p
\end{equation}
except for the speedup of the redistribution part of the program, which is defined with respect to the redistribution time on 2 processes instead of 1 (as there is no redistribution on 1 process).
{\color{black}
A performance reduction in the redistribution phase is visible  on the right in Fig.~\ref{fig:scaling}, consistently with the weak scaling results on the left.
Improved adaptation on many processes and parallel redistribution optimization
} are the focus of ongoing work. Both the weak and strong scaling tests have been performed with the release v1.3.0 of \texttt{ParMmg}\cite{ParMmgSW}.

\begin{table}\caption{Mesh statistics for the \texttt{ParMmg} weak scaling test. Number of vertices $ n_v $ and tetrahedra $ n_e $ in input and output using $ p $ processors.}\label{tab:weakTable}
  \center
  \tabcolsep=0.05cm
  \begin{tabular}{c|c c c c c c c c}
  \hline\noalign{\smallskip}
    $ p $ &
    $ n_v^{in} /p $ & $ n_v^{out}/p $ & $ n_v^{out} / n_v^{in} $ &
    $ n_v^{out} $ &
    $ n_e^{in} /p $ & $ n_e^{out}/p $ & $ n_e^{out} / n_e^{in} $ &
    $ n_e^{out} $ \\
   \noalign{\smallskip}\hline\noalign{\smallskip}
2 & 3625 & 1293637 & 356.81 & 2587274 & 18876 & 7780974 & 412.21 & 15561948 \\
4 & 3467 & 1341637 & 386.88 & 5366549 & 18798 & 8072081 & 429.39 & 32288325 \\
8 & 3346 & 1380055 & 412.38 & 11040444 & 18666 & 8306084 & 444.98 & 66448675 \\
16 & 3264 & 1412516 & 432.66 & 22600269 & 18599 & 8503695 & 457.2 & 136059129 \\
32 & 3190 & 1437267 & 450.46 & 45992552 & 18557 & 8654569 & 466.36 & 276946210 \\
64 & 3214 & 1431186 & 445.2 & 91595935 & 18625 & 8619098 & 462.75 & 551622317 \\
128 & 3215 & 1444674 & 449.31 & 184918370 & 18878 & 8701524 & 460.91 & 1113795077 \\
256 & 3345 & 1468905 & 439.01 & 376039759 & 19705 & 8848464 & 449.03 & 2265206835 \\
512 & 3375 & 1446532 & 428.52 & 740624790 & 19998 & 8714450 & 435.74 & 4461798709 \\
1024 & 3335 & 1449215 & 434.54 & 1483996788 & 19821 & 8731162 & 440.49 & 8940710661 \\
\noalign{\smallskip}\hline\noalign{\smallskip}
  \end{tabular}
\end{table}


\begin{table}\caption{Mesh statistics for the \texttt{ParMmg} strong scaling test. Percentage of edges $ N^{(a,b]} $ whose length in the assigned metrics falls in the interval  $ l_{\mathcal{M}} \in (a,b] $, for each simulation on $ p $ processors.}\label{tab:strongTable}
  \center
  \tabcolsep=0.05cm
  \begin{tabular}{c|c c c c c c c c c}
    \hline\noalign{\smallskip}
    $ p $ &
    $ N^{(0   , 0.3 ]}  $ &
    $ N^{(0.3 , 0.6 ]} $ &
    $ N^{(0.6 , 0.7 ]} $ &
    $ N^{(0.71, 0.9 ]} $ &
    $ N^{(0.9 , 1.3 ]} $ &
    $ N^{(1.3 , 1.41]} $ &
    $ N^{(1.41, 2   ]} $ &
    $ N^{(2   , 5   ]} $ &
    $ N^{> 5         } $ \\
    \noalign{\smallskip}\hline\noalign{\smallskip}
 1 &   1.34 \% &  35.34 \% &  16.89 \% &  25.49 \% &  20.15 \% &   0.63 \% &   0.16 \% & 0  & 0 \\
 2 &   1.14 \% &  34.87 \% &  16.97 \% &  25.79 \% &  20.45 \% &   0.63 \% &   0.15 \% & 0  & 0 \\
 4 &   1.04 \% &  34.20 \% &  17.03 \% &  26.08 \% &  20.85 \% &   0.64 \% &   0.16 \% & 0  & 0 \\
 8 &   0.98 \% &  33.66 \% &  17.06 \% &  26.30 \% &  21.18 \% &   0.66 \% &   0.16 \% & 0  & 0 \\
 16 &   1.07 \% &  32.41 \% &  17.06 \% &  26.78 \% &  21.84 \% &   0.68 \% &   0.16 \% & 0  & 0 \\
 32 &   0.91 \% &  30.78 \% &  17.29 \% &  27.59 \% &  22.57 \% &   0.70 \% &   0.16 \% & 0  & 0 \\
 64 &   0.93 \% &  30.22 \% &  17.20 \% &  27.76 \% &  23.01 \% &   0.72 \% &   0.17 \% & 0  & 0 \\
 128 &   0.88 \% &  29.48 \% &  17.20 \% &  28.06 \% &  23.48 \% &   0.73 \% &   0.17 \% & < 0.01 \% & 0 \\
 256 &   0.70 \% &  27.63 \% &  17.20 \% &  28.84 \% &  24.67 \% &   0.77 \% &   0.18 \% & < 0.01 \% & < 0.01 \%\\
 512 &   0.66 \% &  27.19 \% &  17.14 \% &  28.96 \% &  25.04 \% &   0.80 \% &   0.20 \% & < 0.01 \% & < 0.01 \%\\
 1024 &   0.69 \% &  26.14 \% &  16.91 \% &  29.03 \% &  25.99 \% &   0.92 \% &   0.30 \% &   0.02 \% & < 0.01 \%\\
  \noalign{\smallskip}\hline\noalign{\smallskip}
  \end{tabular}
\end{table}

\renewcommand{\SketchScale}{0.4\textwidth}

\begin{figure}[htbp]
    \centering
     \includegraphics[width=0.42\textwidth,trim=0mm 0mm 0mm 0mm,clip]{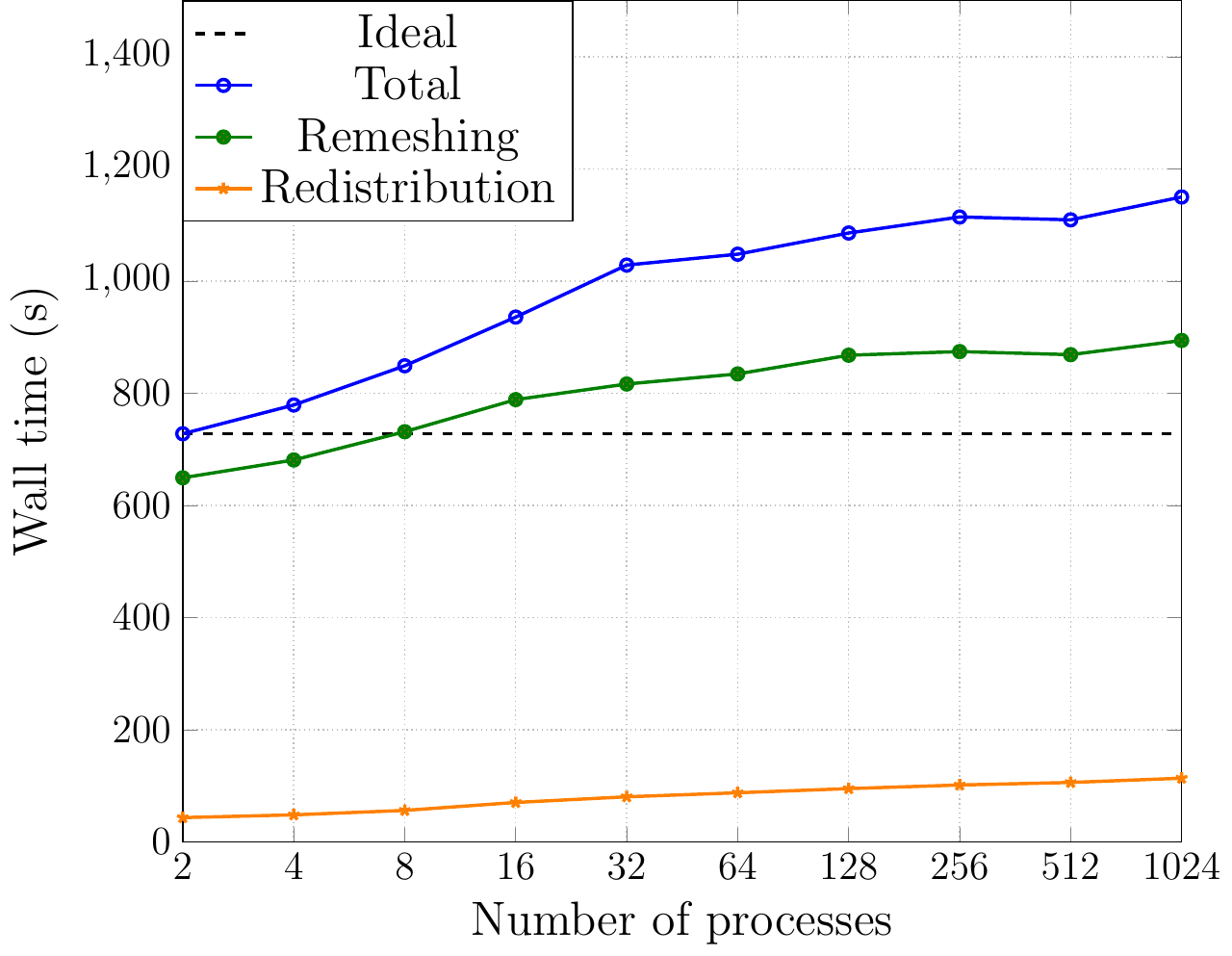} 
     \includegraphics[width=\SketchScale,trim=0mm 0mm 0mm 0mm,clip]{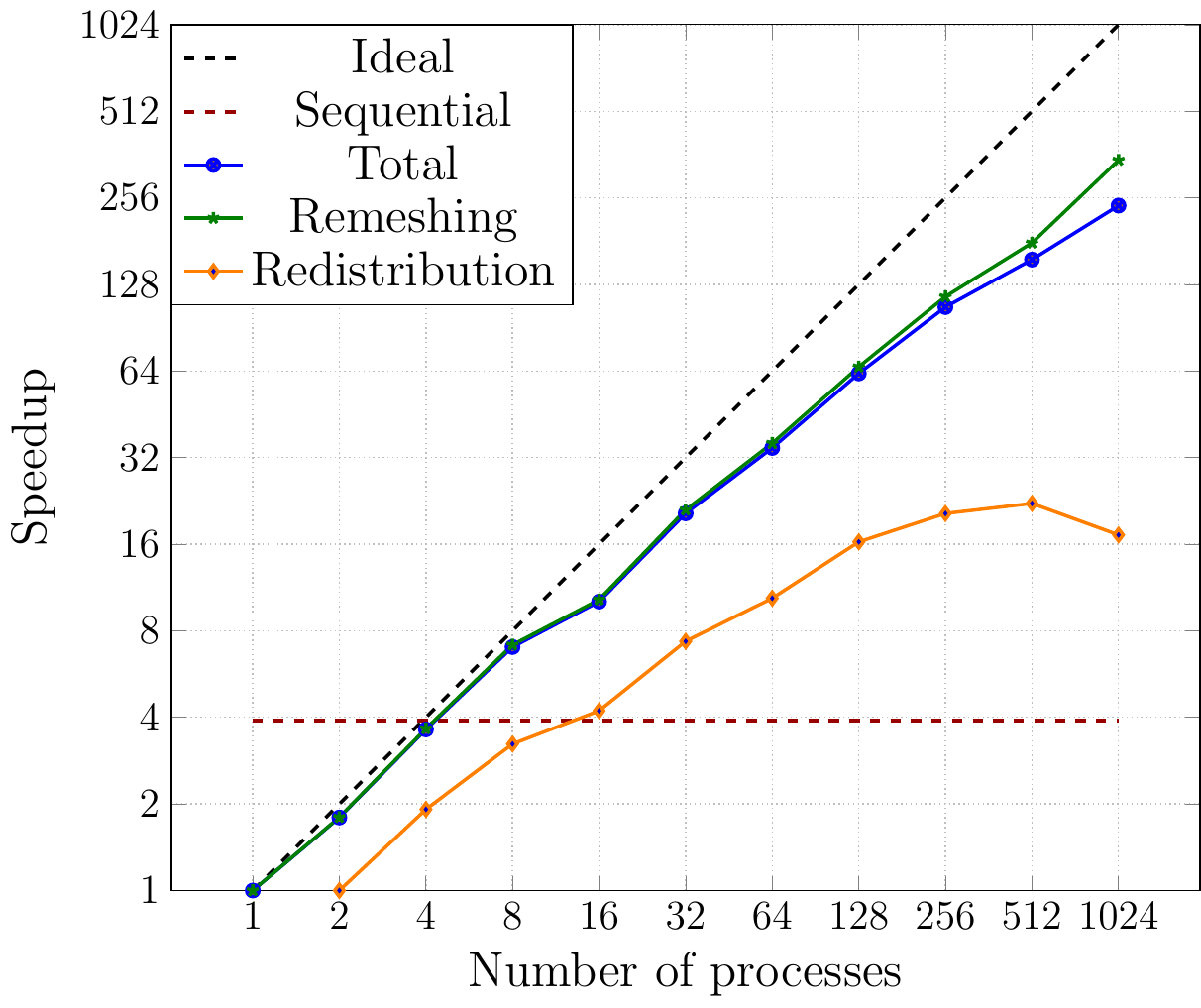} 
    \caption{\texttt{ParMmg} scaling test. \textbf{left} Weak scaling for the uniform refinement of a sphere. \textbf{right} Strong scaling for the adaptation to the double Archimedean spiral in a sphere.}\label{fig:scaling}
\end{figure}

\renewcommand{\SketchScale}{0.4\textwidth}

\begin{figure}[htbp]
    \centering
    \subfloat[Surface adaptation (1 process)]{
      \includegraphics[width=\SketchScale,trim=0mm 8mm 0mm 8mm,clip]{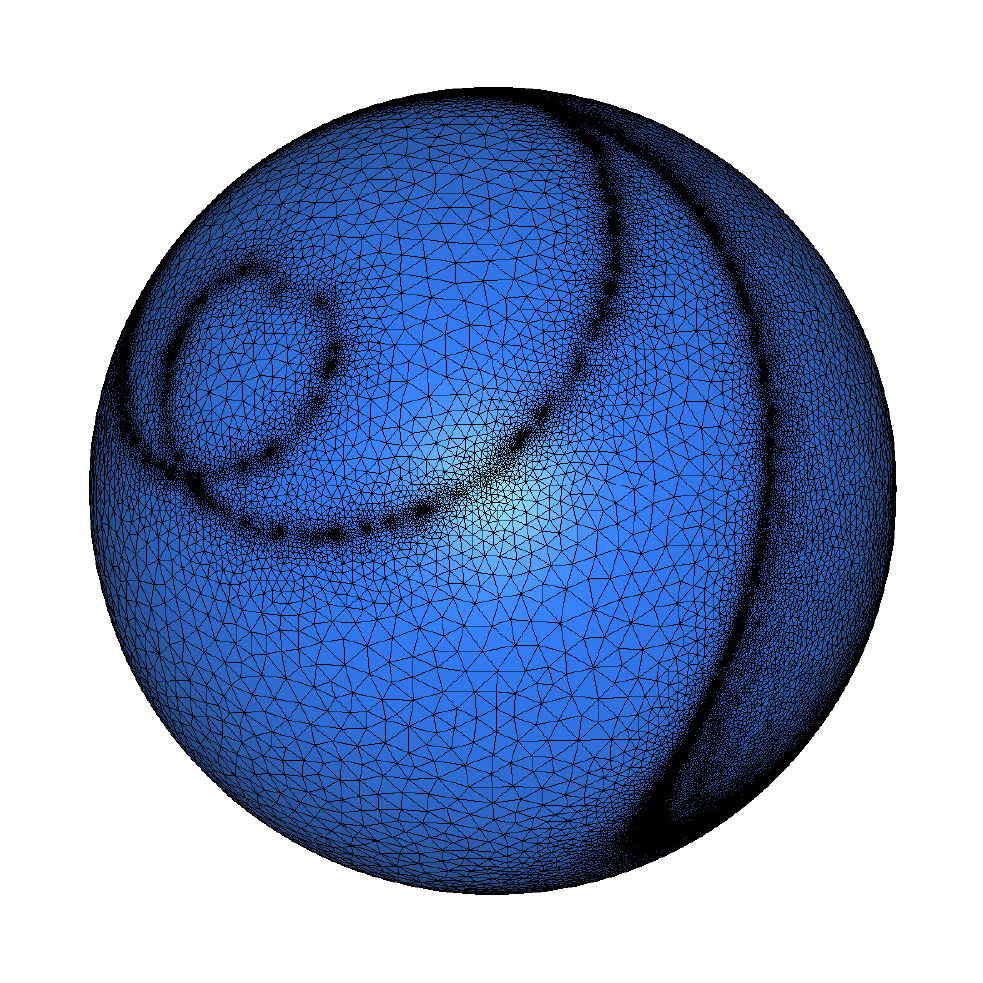}\label{fig:spiral_p1_surf}}
    \subfloat[Surface adaptation (1024 processes)]{
      \includegraphics[width=\SketchScale,trim=0mm 8mm 0mm 8mm,clip]{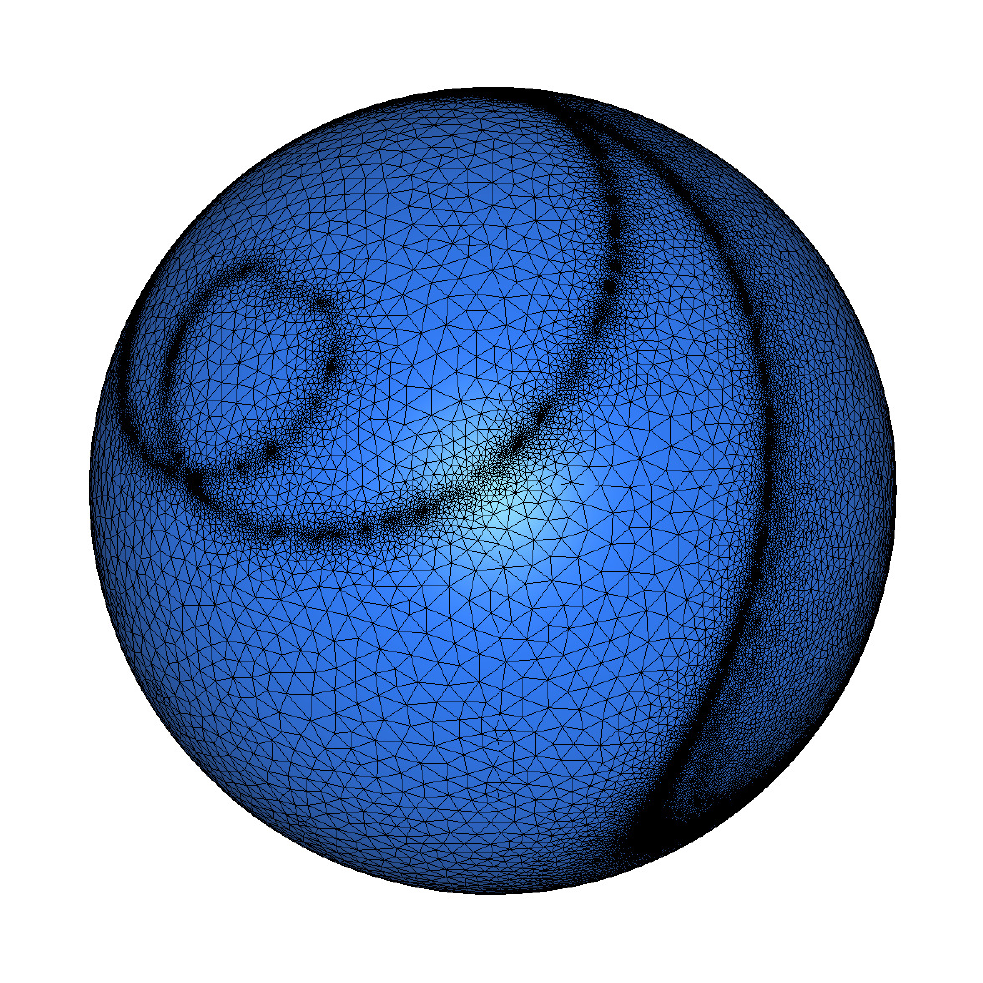}\label{fig:spiral_p1024_surf}}\\
    \subfloat[Volume adaptation (1 process)]{
      \includegraphics[width=\SketchScale,trim=0mm 8mm 0mm 12mm,clip]{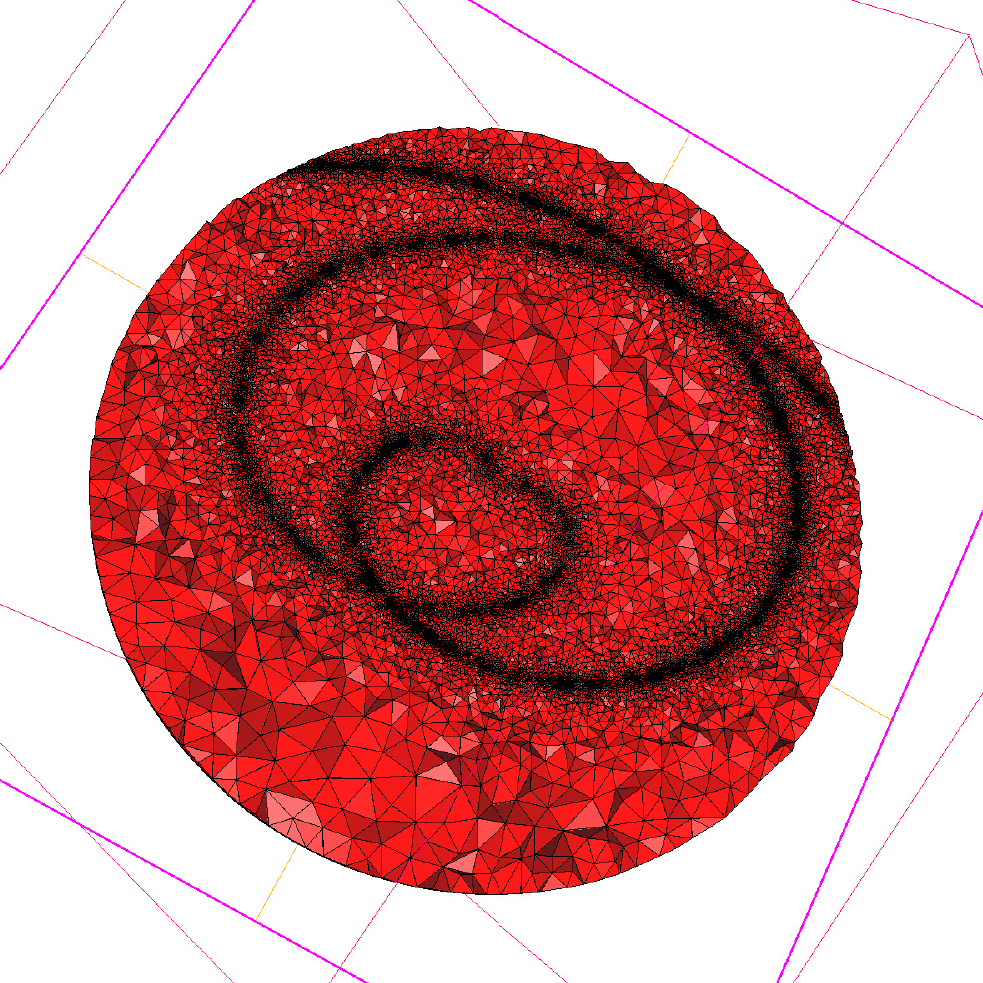}\label{fig:spiral_p1_cut}}
    \subfloat[Volume adaptation (1024 processes)]{
      \includegraphics[width=\SketchScale,trim=0mm 8mm 0mm 12mm,clip]{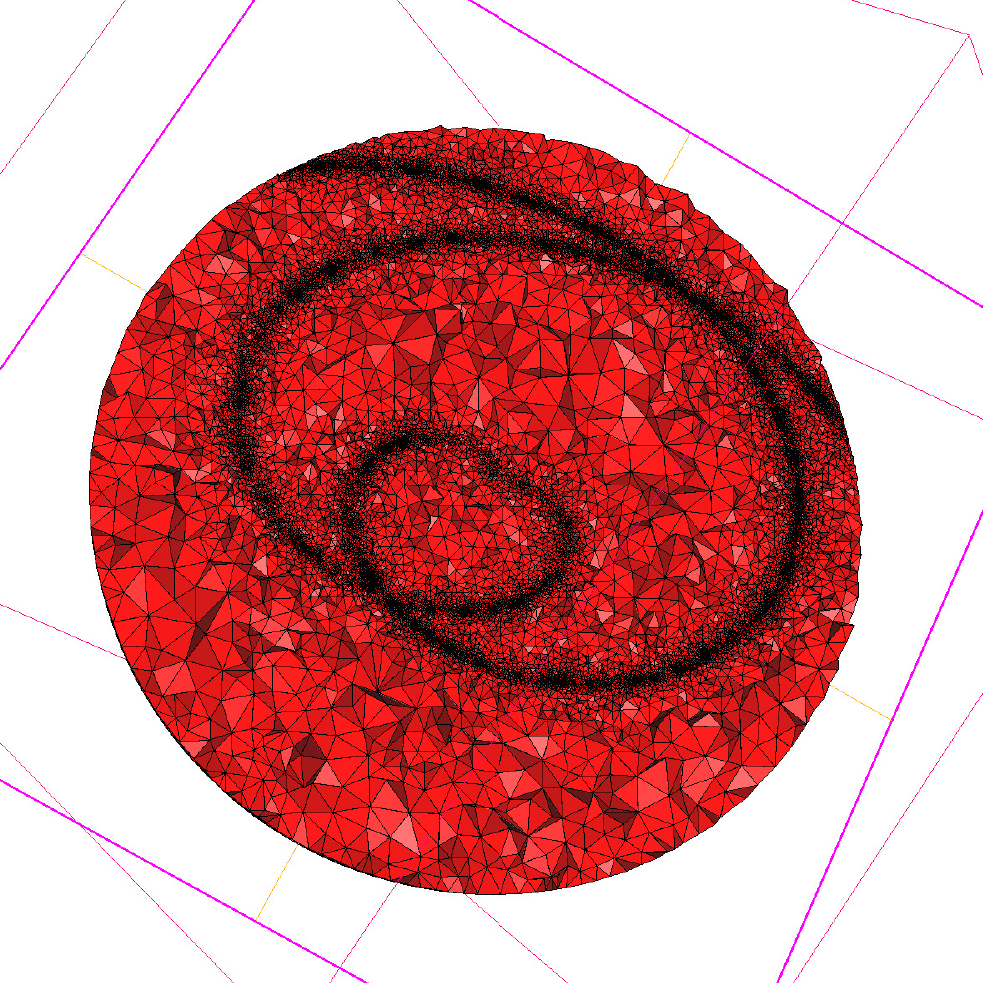}\label{fig:spiral_p1024_cut}}\\
    \caption{Adaptation to the double Archimedean spiral on 1 and 1024 processes.}\label{fig:Archimede}
\end{figure}

%% file: rgeneral.tex
$r-$adaptation techniques are an appealing approach  for unsteady simulations of sharp moving fronts.
Compared to      $h-$refinement,  one of their attractive characteristics is   the relative algorithmic simplicity,
and the  ease of defining conservative remaps due to the inherent continuous nature of the process.

Complete reviews of $r-$adaptation can be found for example in~\cite{Knupp1993,Tang2005,Budd2009,huang-russel2011}.   
In these methods mesh movement is based on parabolic or elliptic differential maps,   often referred to as  mesh PDEs.  
Several interesting formulations allowing PDE based mesh movement have also been proposed for fluid structure interaction
where formulations in which the mesh is treated as an  elastic continuum  are often employed  (see e.g. \cite{Stein2003} and references therein).
Although the problem can be formulated as   a mapping $ \xibf : \Omega_{\xbf} \rightarrow \Omega_{\xibf} $ from the physical domain to a reference one,
for computational efficiency~\cite{Ceniceros2001}  the PDE problem is often  written  for the inverse mapping $ \xbf : \Omega_{\xibf} \rightarrow \Omega_{\xbf} $.
Thus, a solution is sought directly for the mesh nodes position $ \xbf $, and the mesh PDE is discretized  in the (fixed) reference domain $ \Omega_{\xibf} $. 
This framework allows a direct analogy between the mesh PDE and  Lagrangian or Arbitrary Lagrangian Eulerian continuum mechanics, and its associated computational methods.

Introducing  the displacement $ \deltabf = \xbf -\xibf $, the mesh PDE problem can be written as 
\begin{equation}
\label{eq:mmpde0}
\tau \partial_t\deltabf - \grad_{\xibf} \cdot \sigmabf ( \deltabf ) =  \pmb{F}(\xibf,\deltabf)
\end{equation}
where $ \sigmabf ( \deltabf )$ can be seen as a Cauchy stress tensor, and the volume force $\pmb{F}(\xibf,\deltabf)$  
depends on smoothness estimator.  \mario{Problem \eqref{eq:mmpde0} should be  supplemented by  appropriate  boundary conditions.
To  limit the focus of this contribution we leave this aspect out of the discussion.
The interested reader can refer to  \cite{CIRROTTOLA2021110177} 
for a general formulation involving Dirichlet, and slip  conditions along curved boundaries.}

 The parameter/matrix $\tau$  is a   dimensionless relaxation time \mario{allowing} 
 to  adjust the  time scale of the mesh movement  w.r.t. the  physics. This  parabolic  regularisation/relaxation is also 
a good model for  elliptic solvers  coupled with some  truncated iterative process,  often used  in practice especially when  
explicit time stepping is used in the flow solver  \cite{Ceniceros2001,Tang2003,Chen2008,nouveau:tel-01500093,arpaia:hal-01372496,Arpaia2020}.
\mario{In particular, in our case the mesh PDEs are discretized with  standard $\mathbb{P}^1$ finite elements on the reference/initial  mesh,
and    new nodal  positions are obtained via  a small number of  nodal  Jacobi iterations that can be written as 
\begin{equation}\label{eq:Jacobi0}
\deltabf_i^{[k+1]} = \deltabf_i^{[k]} - ( K_{ii}^{[k]} )^{-1}  \sum\limits_{j \in\mathcal{B}_i} K_{ij}^{[k]}\deltabf_j^{[k]}  +  \pmb{F}_i^{[k]}
\end{equation}
with $K$ denoting the stiffness matrix, and $K_{ii}$ being a block or a diagonal matrix depending on the closure for  $ \sigmabf ( \deltabf )$.}\

%% file: rresults2D.tex
\subsection{Embedded geometries and  distance function: Laplace vs elasticity}
\label{sec:rresults}

\mario{The closure laws $ \sigmabf ( \deltabf )$ and $  \pmb{F}(\xibf,\deltabf)$ are an essential part of the method.
There exist a variety of possibilities  that go from 
simple linear Laplace smoothing to nonlinear gradient flow formulations with embedded non-negativity conditions for the Jacobian of the resulting map  \cite{MASUD200777,Budd2009,Huang2001,Huang2015}.}
A key role   is played by the definition of a monitor function/matrix, which is in general a function of the flow solution in the physical space, and plays the role of the metric field in  $h-$adaptation.
In the simplest  setting, if adaptation is performed w.r.t. the field $u$, one starts  by defining a scalar function $\omega=\omega(u(\xbf))$. 
\mario{A somewhat classical definition is} 
\begin{equation}
\omega = ( 1 + \alpha u + \beta  \|\grad_{\xbf}  u\| + \gamma  \|\mathcal{H}_{\xbf}(u)\|)^p
\end{equation}
In this work we consider two cases: 
\begin{enumerate}
\item the Laplacian mesh movement  corresponding  to the closure  $ \sigmabf ( \deltabf )= \omega(\xbf) \grad_{\xibf} \deltabf$ 
\item the elastic closure  $ \sigmabf ( \deltabf )= 2\mu \epsilonbf(\deltabf ) + \lambda \,\text{tr}( \epsilonbf(\deltabf )) $,
with   $(\mu,\lambda)$  constant Lam\'e coefficients, and   $\epsilonbf(\deltabf ) $ the (symmetric) deformation tensor. 
\end{enumerate}
Note that for the Laplacian mesh movement, the classical elliptic equation 
\begin{equation}
\grad_{\xibf} \cdot \left( \omega(\xbf) \grad_{\xibf} \xbf \right) = \zerobf \qquad \text{in}\ \Omega_{\xibf}
\end{equation}
is  obtained for   $ \pmb{F}(\xibf,\deltabf)=\grad_{\xibf} \cdot( \omega(\xbf) \mathbf{Id})$.  With  elasticity,  usually used to propagate 
an imposed boundary  displacement, it is less clear what form the forcing should take. Note that these matrices are constant
throughout the iterations when using linear elasticity.\\

One of our first aims is to enhance as much as possible the  resolution of     surfaces implicitly defined by a level-set function $\Phi$.
We assume to be able to work just with the Dirichlet and natural boundary conditions in Eq.~\ref{eq:mmpde0}.
We use this as an exercise to compare the formulations presented in the previous section.
Although simple, cost-wise the Laplacian movement presents  a considerable non-linearity due to the dependence on the  
final configuration both of the  stress definition, and in the force. As an alternative, we have tried to use linear elasticity as a mean to relax 
the strain locally generated  by  an applied volume force. To this end, we have used the  same force definition for all systems, namely $\pmb{F}(\xibf,\deltabf)=\grad_{\xibf} \cdot( \omega(\xbf) \mathbf{Id})$.

The definition of the monitor function can be  set to mimic a mollified   characteristic function corresponding to the domains  $\Phi >0$ (or $\Phi < 0$).
A possible choice is 
\begin{equation}
\label{eq:omega0}
\omega_{\Phi}^{\text{GB}}=\sqrt{\alpha_0 +\alpha_{\Phi} e^{-\beta_{\Phi} \Phi^2}}
\end{equation}
with GB standing for gradient based. The constant term $\alpha_0$ is   defined to offset the  stiffness in different regions of the domain. A first experiment consists in comparing Laplace and elasticity based
deformation with a similar definition of the sensor. In the experiment we present, the reference domain is meshed with a relatively uniform triangulation, and the Newton-Jacobi update
of Eq.~\ref{eq:Jacobi0} is repeated for a few thousands iterations or   until  three orders of magnitude convergence  are  achieved.

\begin{figure}
\begin{center}
      \includegraphics[angle=90,height=0.33\textwidth]{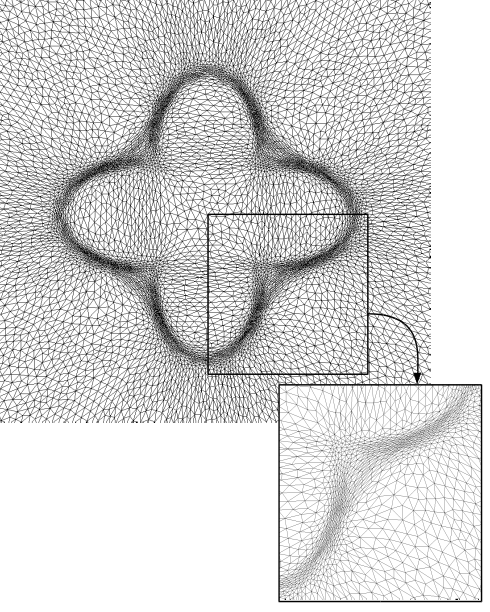}
      \includegraphics[angle=90,height=0.33\textwidth, trim=10 0 0 10,clip]{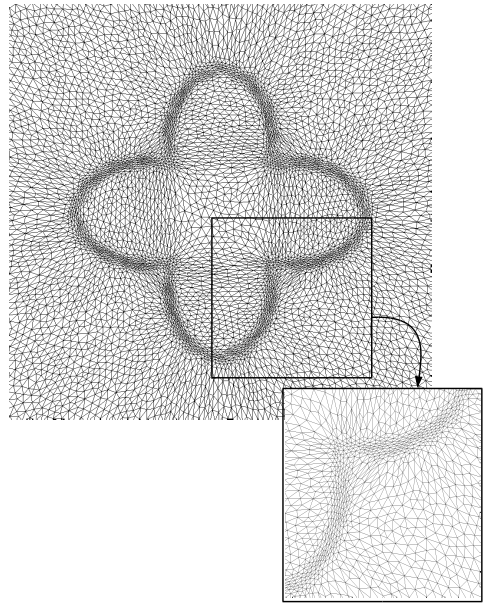}
\caption{\label{fig:flower}Embedded geometries. 4-Lobed Flower adaptation. \textbf{left} Elasticity based movement. \textbf{right} Laplacian based movement.}
\end{center}
\end{figure}

\begin{figure}
\centering
\includegraphics[width=0.3\textwidth,trim=1mm 1mm 1mm 1mm,clip]{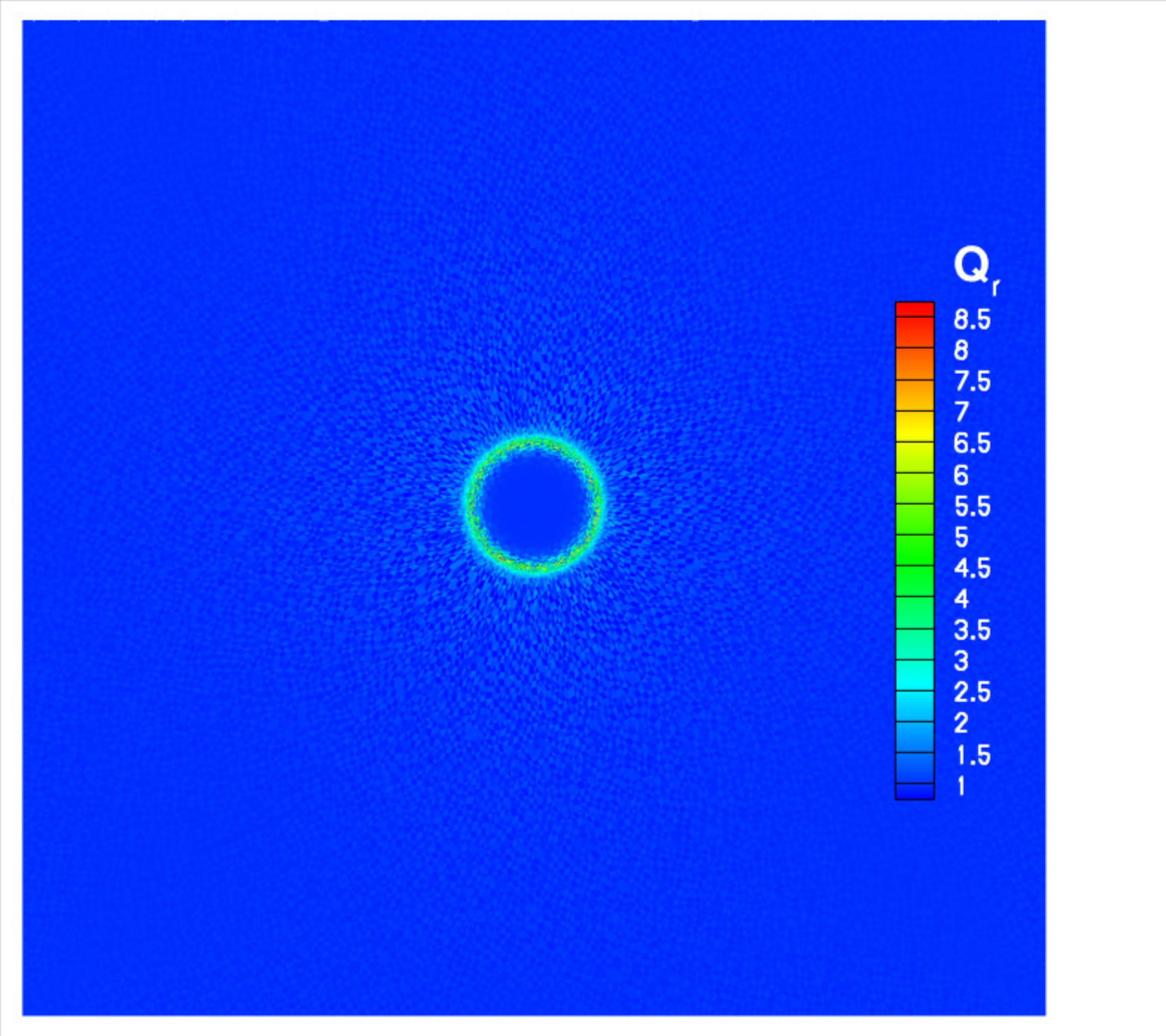}\hspace{0.4cm}
\includegraphics[width=0.3\textwidth,trim=1mm 1mm 1mm 1mm,clip]{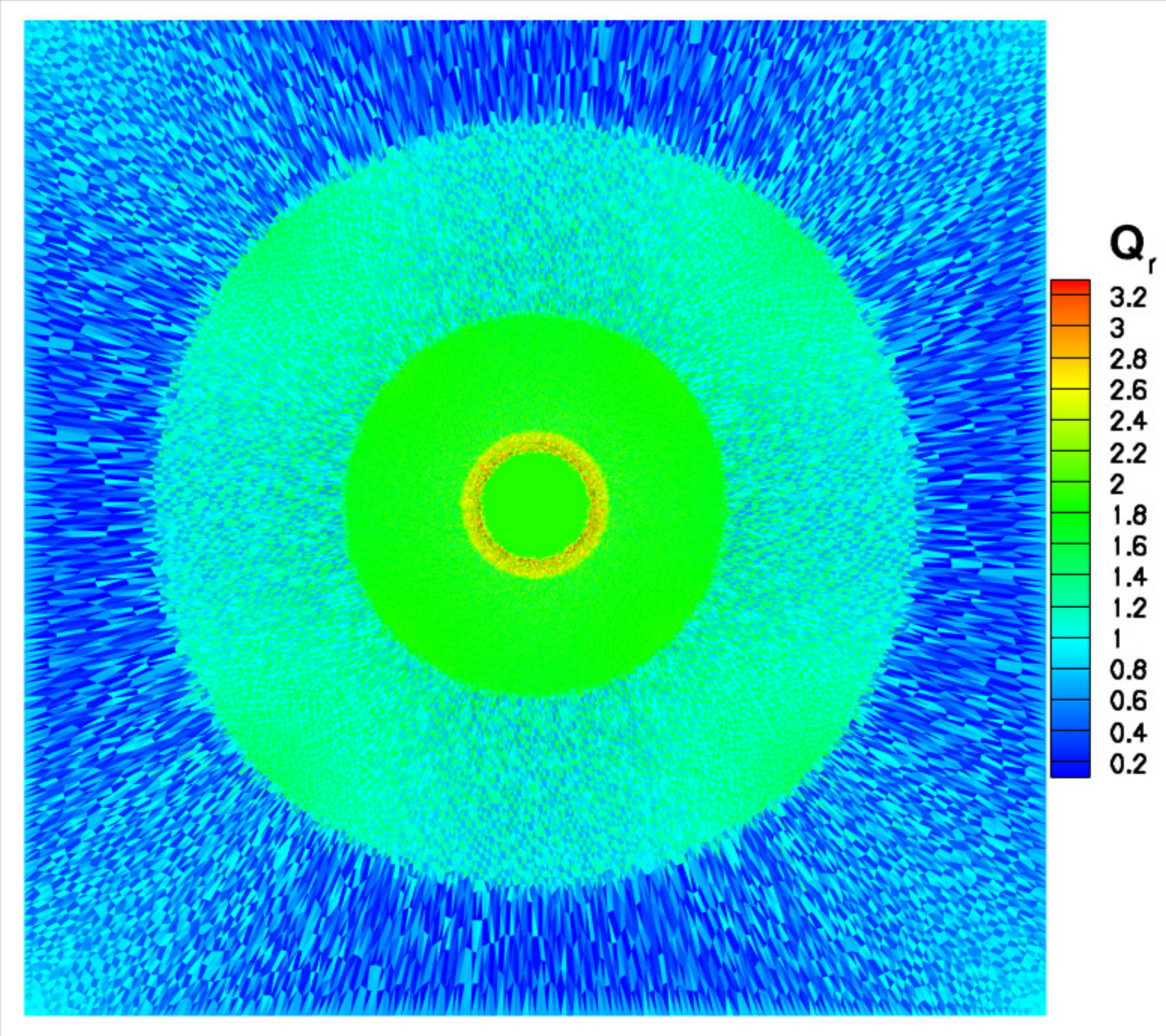}\\[3pt]
\includegraphics[width=0.3\textwidth,trim=1mm 1mm 1mm 1mm,clip]{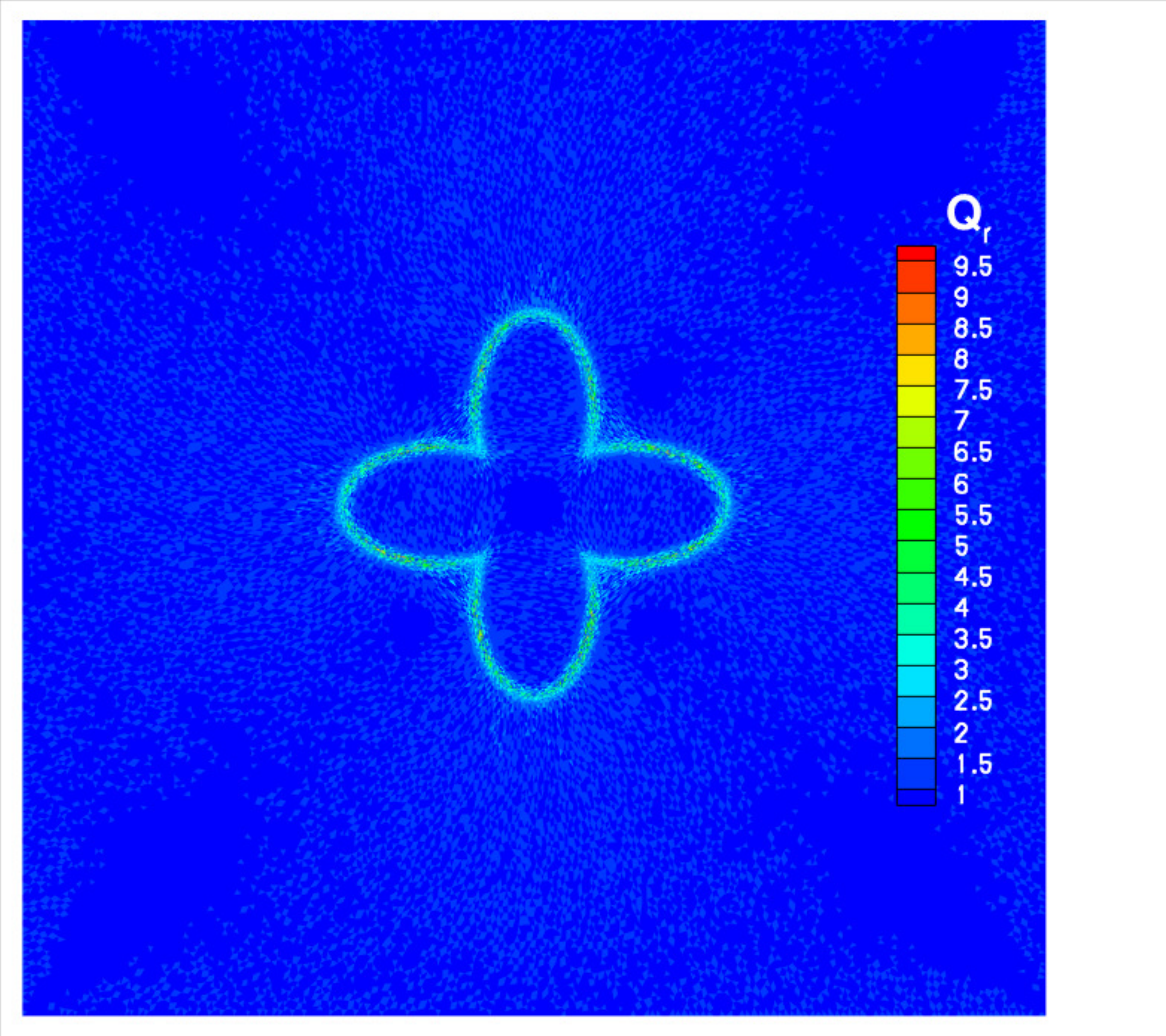}\hspace{0.4cm}
\includegraphics[width=0.3\textwidth,trim=1mm 1mm 1mm 1mm,clip]{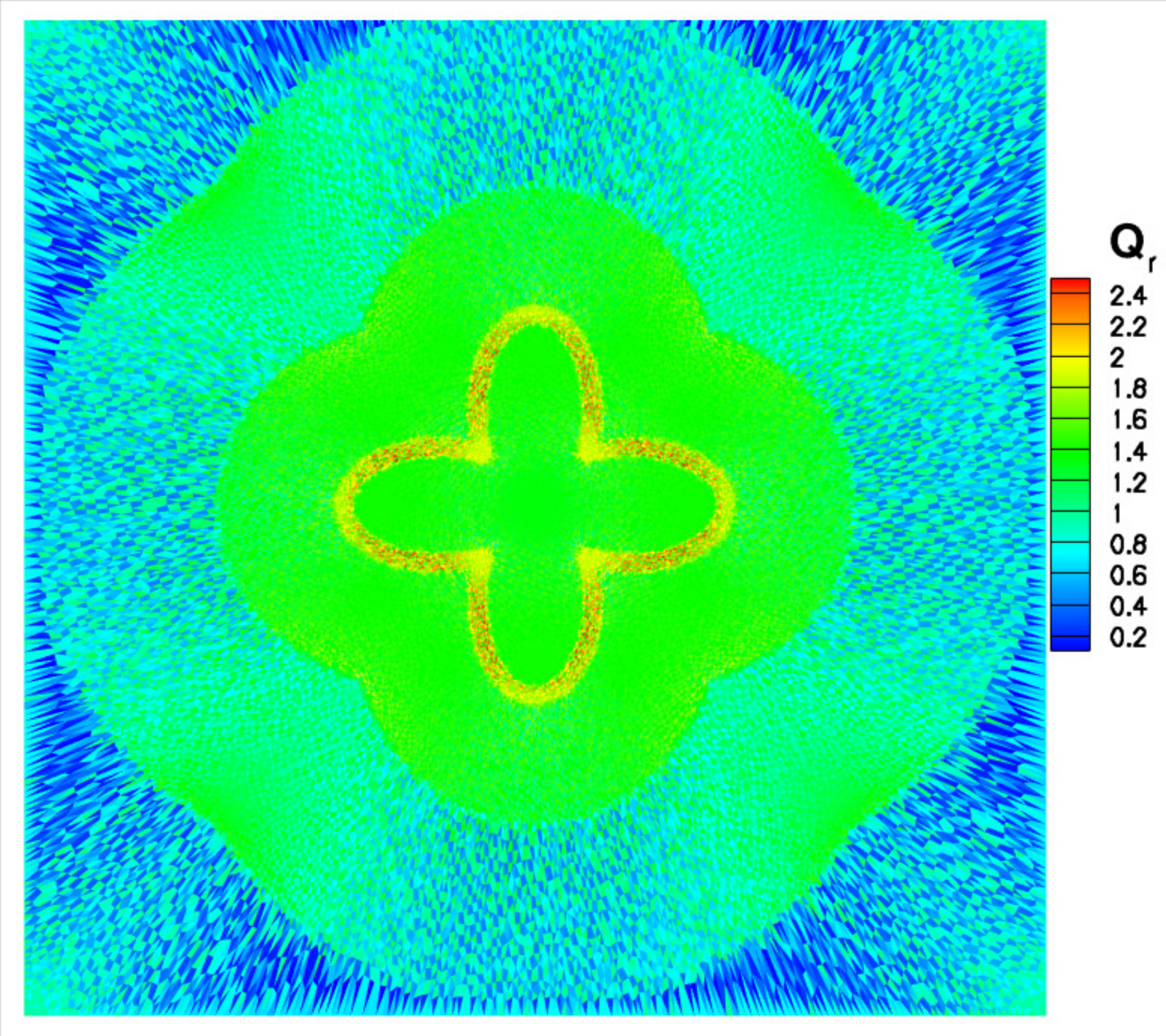}
    \caption{Embedded geometries. Anisotropic compression, $\Qcal_r$, of the adapted mesh for a Circle and a 4-Lobed Flower. \textbf{left} Elasticity. \textbf{right} Laplacian.}\label{fig:compression_embedded}  
\end{figure}

We use as a target the level-set corresponding to an implicit curve resembling a 2D four-lobed flower. The values of the constants in Eq.~\ref{eq:omega0} are slightly tuned to improve the resolution of the curve
but have  similar orders of magnitude in the two approaches:   $(\alpha_{\Phi}, \beta_{\Phi}) = (40,300)$ for the Laplace deformation, and  $(\alpha_{\Phi}, \beta_{\Phi}) = (30,200)$ for elasticity.
The offset is  chose as  $\alpha_0=1$ for the Laplace deformation, while for elasticity we have set $\alpha_0=2.5 |\kappa|$  with $\kappa$ an approximation of the curvature of the level-set contours. 
The latter is introduced in an attempt to improve the capabilities of capturing singularities of the linear elastic deformation.  The two approaches are compared in Fig.~\ref{fig:flower}.
While allowing more anisotropy  in regions with smaller curvature, the elasticity based deformation  tends to provide  highly stretched meshes  in vicinity of singularities as in the corners.
{
  \color{black}
  The elastic approach tends to produce more anisotropic meshes near the level-set, while keeping good isotropy far from it.
  The Laplacian approach, instead, tends to spread streched elements far from the level-set.
  This behaviour is highlighted in Fig.~\ref{fig:compression_embedded}, where the ratio
  \begin{equation}\label{eq:compdef}
    \Qcal_r = \frac{r_{ref}}{r_{adap}}
  \end{equation}
  between the radius of the circle inscribed to the elements on the reference and adapted mesh is depicted for the meshes adapted to a flower level-set and a circular level-set. The statistics shown in Table~\ref{tab:flowerStats} and Table~\ref{tab:circleStats} also show these trends by reporting, for the elements in a narrow band $ | LS | < 10^{-2} $ near the zero level-set,  the number of elements, the minimum, maximum and average edge sizes and the 2D isotropic quality $ \Qcal_{iso} $ of each element $ K $
  \begin{equation}
    \Qcal_{iso} = \frac{\beta \sum_{j = 1}^3 L_j^2}{S_K}
  \end{equation}
  Here the isotropic quality is purposely used to highlight the increased capability of the elastic approach to produce anisotropic elements near the level-set.} This may be a limiting factor in more complex situations as e.g. when combining level-set with solution based adaptation.

\begin{table}\caption{Circle adaptation mesh statistics.}\label{tab:circleStats}
  \center
  \tabcolsep=0.05cm
  \begin{tabular}{c c c c c c c c c}
    \hline\noalign{\smallskip}
                 &  \# elements & Min. $h$ & Max. $h$ & Av. $h$ & Min. $\Qcal_{iso}$ & Max. $\Qcal_{iso}$ & Av. $\Qcal_{iso}$ \\
    \noalign{\smallskip}\hline\noalign{\smallskip}
    Initial mesh &  203 & $ 9.19 \times 10^{-3} $ & $ 2.60 \times 10^{-2} $ & $ 1.58 \times 10^{-2} $ & 1 & 1.33 & 1.05 \\
    Laplacian    & 1758 & $ 3.72 \times 10^{-3} $ & $ 0.15                $ & $ 6.57 \times 10^{-3} $ & 1 & 1.71 & 1.12 \\
    Elasticity   & 3379 & $ 1.71 \times 10^{-3} $ & $ 3.79 \times 10^{-2} $ & $ 6.11 \times 10^{-3} $ & 1 & 6.50 & 2.42 \\
    \noalign{\smallskip}\hline\noalign{\smallskip}
  \end{tabular}
\end{table}

\begin{table}\caption{4-Lobed Flower adaptation mesh statistics.}\label{tab:flowerStats}
  \center
  \tabcolsep=0.05cm
  \begin{tabular}{c c c c c c c c c}
    \hline\noalign{\smallskip}
                 &  \# elements & Min. $h$ & Max. $h$ & Av. $h$ & Min. $\Qcal_{iso}$ & Max. $\Qcal_{iso}$ & Av. $\Qcal_{iso}$ \\
    \noalign{\smallskip}\hline\noalign{\smallskip}
    Initial mesh &  655 & $ 9.19 \times 10^{-3} $ & $ 2.60 \times 10^{-2} $ & $ 1.58 \times 10^{-2} $ & 1 & 1.33 & 1.05 \\
    Laplacian    & 3195 & $ 4.67 \times 10^{-3} $ & $ 0.15                $ & $ 8.74 \times 10^{-3} $ & 1 & 1.85 & 1.15 \\
    Elasticity   & 7793 & $ 1.49 \times 10^{-3} $ & $ 4.05 \times 10^{-2} $ & $ 7.1  \times 10^{-3} $ & 1 & 7.90 & 2.52 \\
    \noalign{\smallskip}\hline\noalign{\smallskip}
  \end{tabular}
\end{table}

To  obtain a high  resolution of the zero level-set without impoverishing too much the surrounding regions,  
the Laplacian movement could be combined with some 
 way to directly control the mesh density  at a given distance from the embedded surface. A possible choice is a piecewise constant definition of the monitor function:
\begin{equation}\label{eq:omega1}
\omega_{\Phi}^{\text{PC}}=\left\{
\begin{array}{lll}
\omega_1 \quad&\text{if }\;\; |\Phi| \le \Phi_1\\[5pt]
\omega_j \quad&\text{if }\;\;    \Phi_{j-1}  < |\Phi|\le\Phi_{j}\;\;\forall j\in{2,N_{\Phi}-1}\\[5pt]
\omega_{N_{\Phi}}\quad&\text{if }\;\;     \Phi_{N_{\Phi}-1}  < |\Phi| 
\end{array}
\right.
\end{equation}
The number of regions $N_{\Phi}$ is usually taken at least equal to three. This definition allows to set up areas of increasing stiffness, thus with higher density of nodes, 
as one approaches  the zero level-set. The last jump around the usually small value $\Phi_0$ provides  the refinement sought  around the embedded surface.
In this area  \eqref{eq:omega0}  could also be used, but a higher constant value works fairly well.

\begin{figure}
      \includegraphics[width=0.28\textwidth, trim=220 240 220 240,clip]{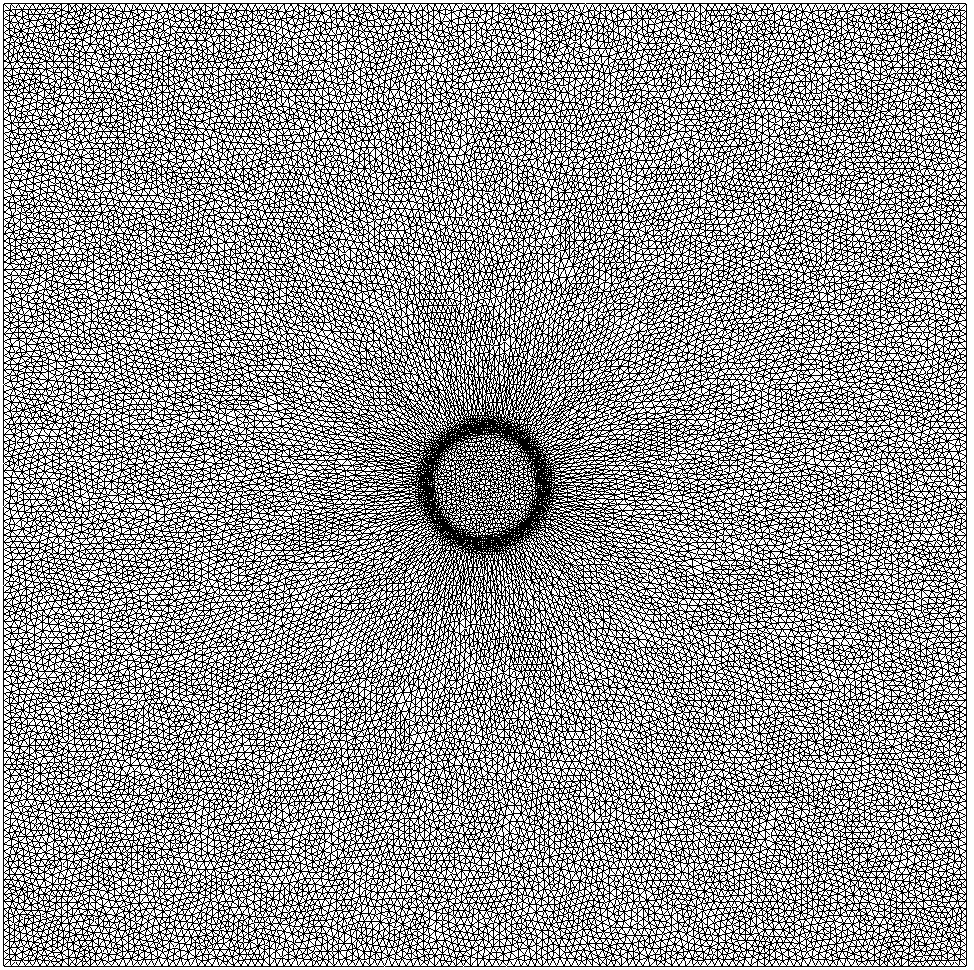}
      \includegraphics[width=0.28\textwidth, trim=220 240 220 240,clip]{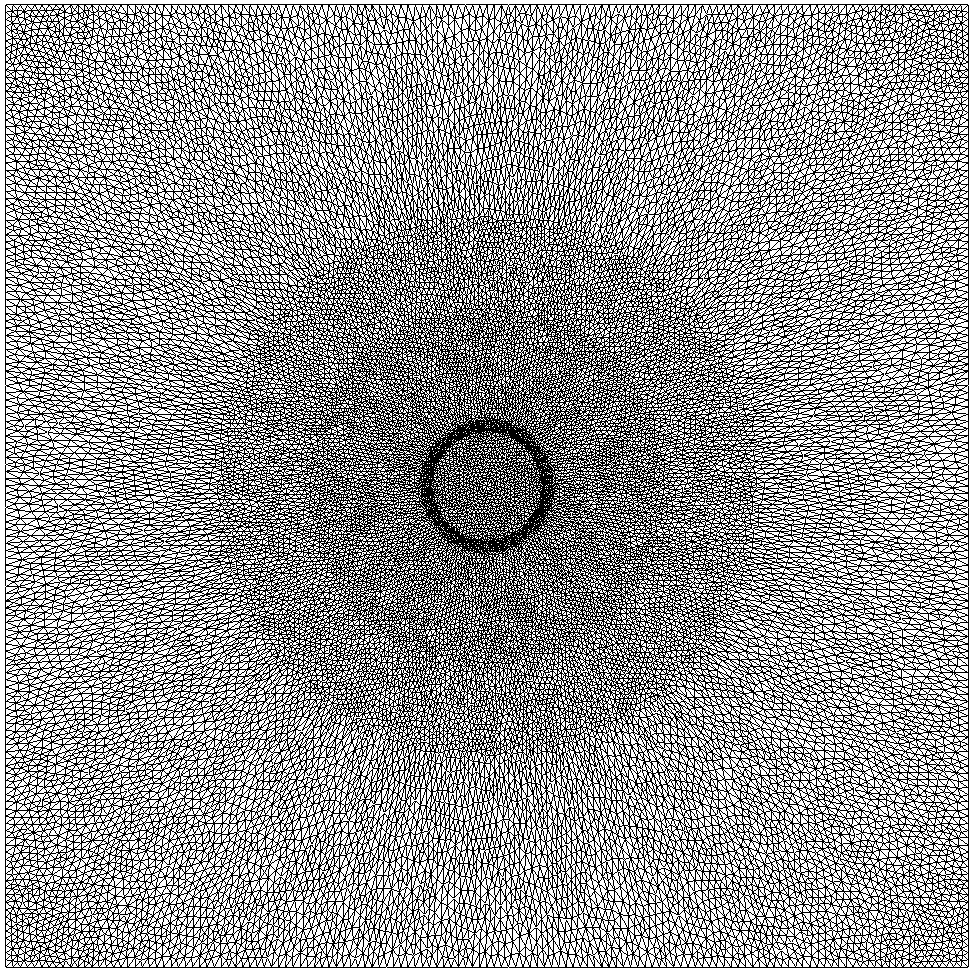}
      \includegraphics[width=0.415\textwidth]{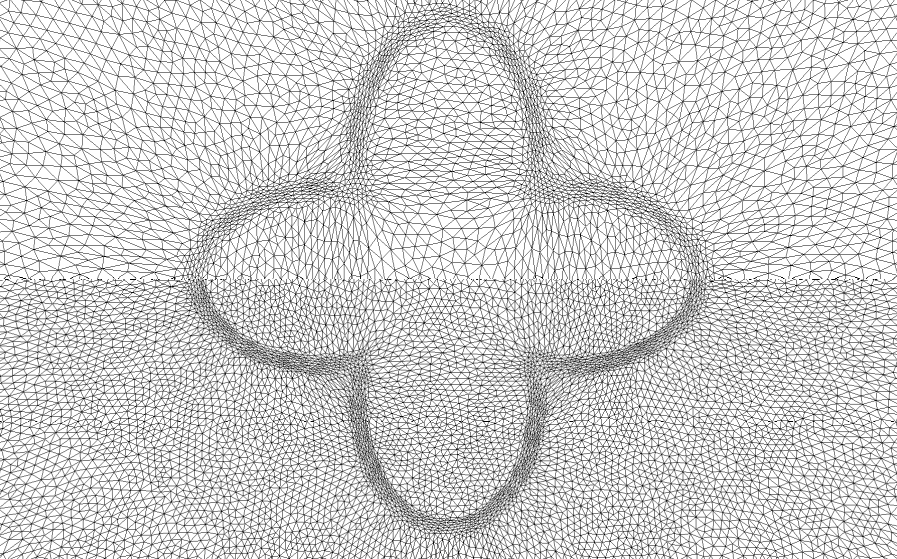}\\
\caption{\label{fig:flower1}Embedded geometries. Gradient based (GB) versus piecewise constant (PC) monitor function. \textbf{left} Circle: Laplacian GB. \textbf{center} Circle: Laplacian PC. \textbf{right} 4-Lobed Flower: Laplacian GB (top) -  Laplacian PC (bottom).}
\end{figure}

To give an example we consider  a circular interface,   and we compare on the left and central pictures of Fig.~\ref{fig:flower1}
the gradient based (or GB) definition with  a 4 layers PC with : $(\Phi_1,\Phi_2,\Phi_3)=
(0.05,1,1.75)$, and $(\omega_1,\omega_2,\omega_3)=
(225,90,70,20)$. We see that this definition allows to  obtain a denser area  of uniform mesh size closer to the zero level-set,
thus virtually leaving out nodes for resolving other features, and allows a higher  resolution in correspondence of the  embedded circle.
As a final comparison, in the rightmost  picture in Fig.~\ref{fig:flower1}  we show the two approaches compared for the 4-lobed flower.
This simple piece-wise constant definition of the monitor, combined with Laplacian deformation provides  a good degree of control on the 
node density.

In the next paragraph we focus on applications based on this approach. The reader can refer to the PhD thesis~\cite{nouveau:tel-01500093} for more comparisons
and applications involving the elastic deformation using the above definitions, and to~\cite{Huang2001,Budd2009,huang-russel2011} for  the   derivation of mesh PDEs 
which attempt to  impose a  given metric field without topology change.


\subsection{Immersed boundaries and moving bodies}

We consider the combination of  level-set adaptation,  immersed boundary methods, and solution based adaptation.
The PDE setting is  similar to that of Sect.~\ref{sec:ibmresults}, but the numerical results use a  flow solver based on
a stabilized residual based method    in an Arbitrary-Lagrangian-Eulerian  (ALE) setting~\cite{nouveau:hal-01348902,nouveau:tel-01500093,AR:17}. 
We focus on forced motion of a solid, whose boundary  is  implicitly defined via the zero level-set of a signed distance function.
To cope with the time dependent nature of the problem in this case  we   use the explicit knowledge of the exact  geometry,
and recompute the distance using the method proposed in \cite{mshdist-dapogny}. As before, we assume to be able to work just with the Dirichlet and natural boundary conditions in Eq.~\ref{eq:mmpde0}. \\

As in Sect.~\ref{sec:ibmresults}, we wish to combine  level-set and solution based adaptation.  To this end, we introduce the smoothness monitor for the solution
\begin{equation}\label{eq:omega2}
\omega_{u} =  \sqrt{1+\alpha_u  \widehat{\nabla u}^2 }
\end{equation}
where $\widehat{\nabla u}$ is a capped gradient norm (see e.g.   \cite{Tang2003, Chen2008})
\begin{equation}\label{eq:nablac-cap}
\widehat{\nabla u} = \min(1 , \dfrac{ \|\nabla_{\xbf} u \|}{\beta_u \| \nabla_{\xbf} u\|_{\text{ref}} })
\end{equation}
with the reference value $ \| \nabla u\|_{\text{ref}} $ usually chosen as  the maximum over the domain.  Note that  the gradient is computed w.r.t. the actual coordinates $\xbf$. 
The combined level-set/solution monitor function is finally defined as 
\begin{equation}\label{eq:omega3}
\omega =\left\{
\begin{array}{ll}
\omega_{\Phi} \quad&\text{if }\;\; |\Phi| \le \epsilon \\[5pt]
\max(\omega_{\Phi},\omega_{u})\quad&\text{if }\;\; |\Phi|  >\epsilon  
\end{array}
\right.
\end{equation} 
The initial mesh is adapted to the embedded surface and initial solution, solving the Newton-Jacobi iterations until a relative convergence of a few orders of magnitude,
 as discussed in the previous  section. However, during the ALE time advancement,  only a few mesh iterations are performed. For the example shown we use 
 20 mesh iterations, with the last known mesh as initial guess. This gives 
 a negligible overhead compared to the implicit Crank-Nicolson ALE time stepping used for the flow.

\begin{figure}
\centering
\includegraphics[height = 3cm]{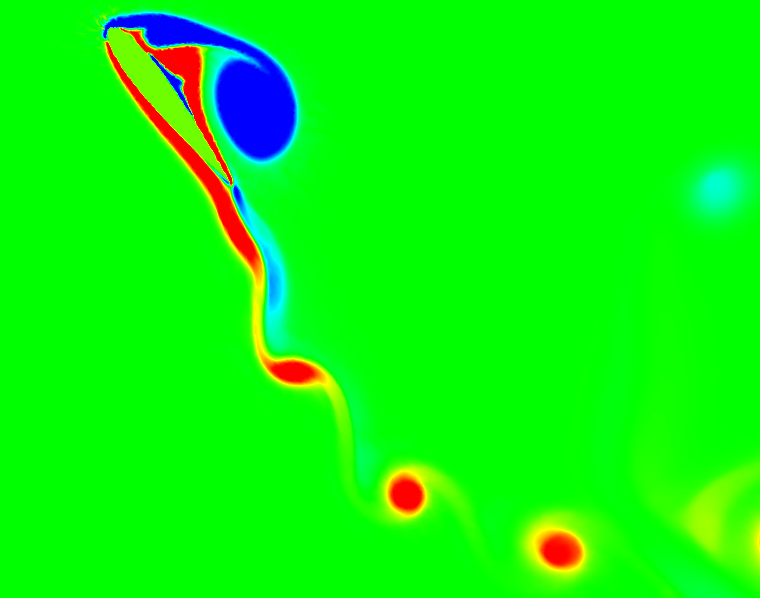}\hspace{0.05cm}
\includegraphics[height = 3cm]{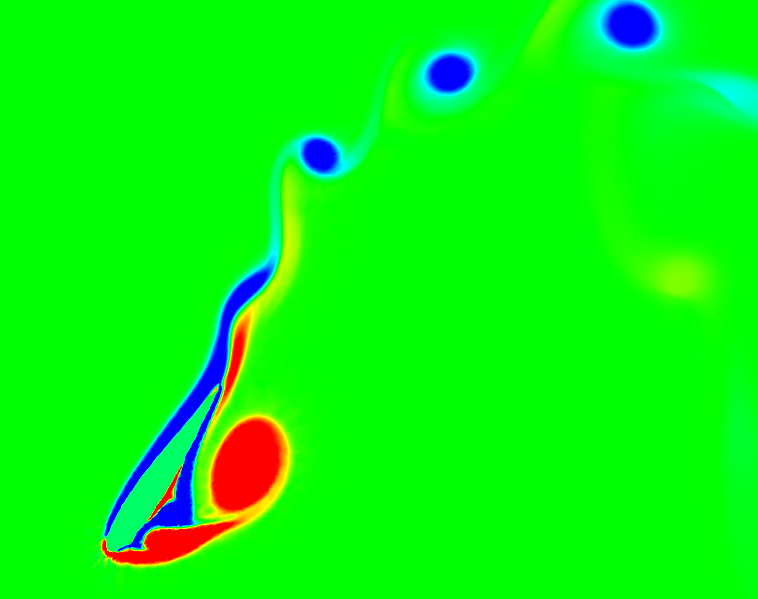}\hspace{0.1cm}
\includegraphics[height = 3cm]{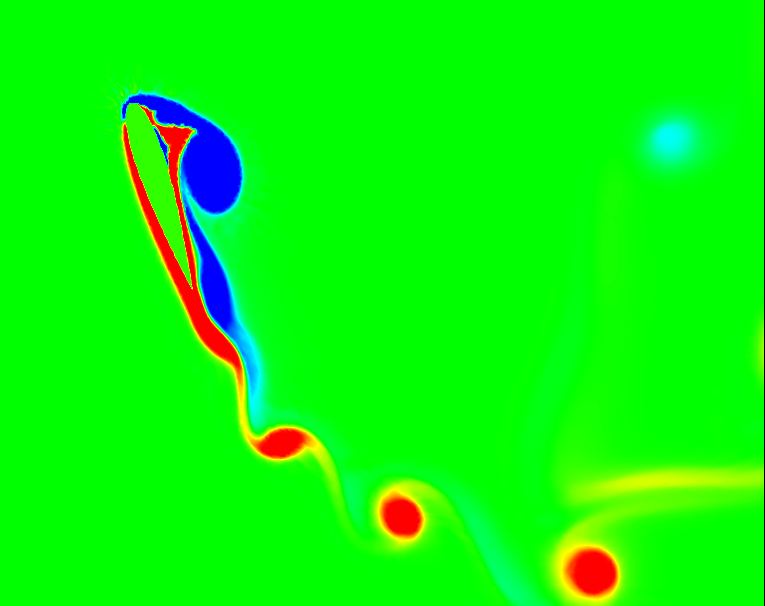}\\
\vspace{0.25cm}

\includegraphics[height = 3cm]{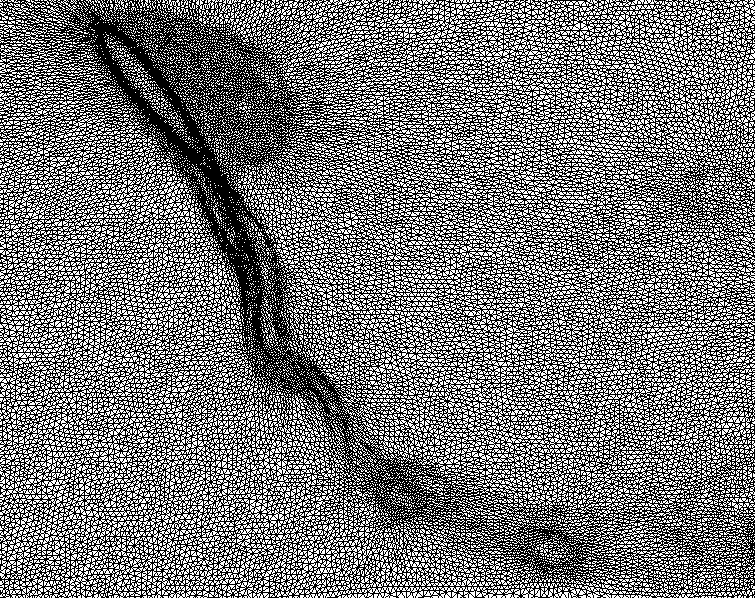}\hspace{0.05cm}
\includegraphics[height = 3cm]{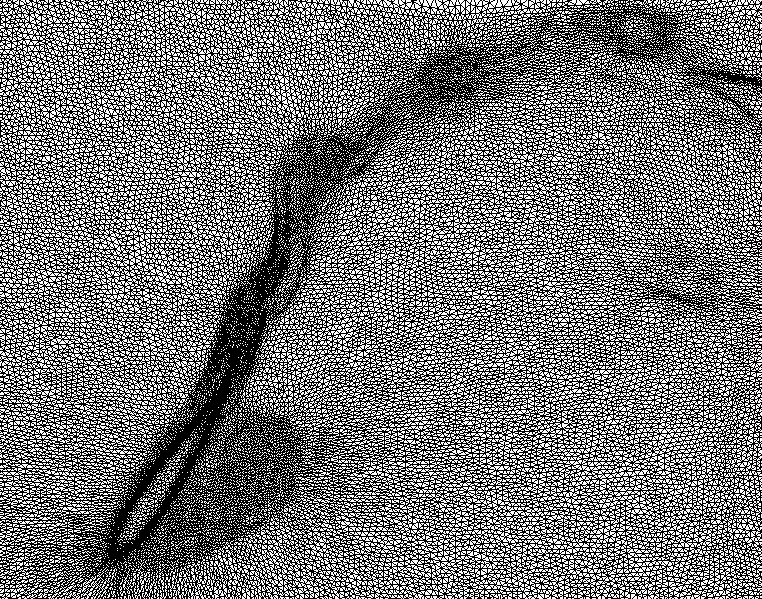}\hspace{0.1cm}
\includegraphics[height = 3cm]{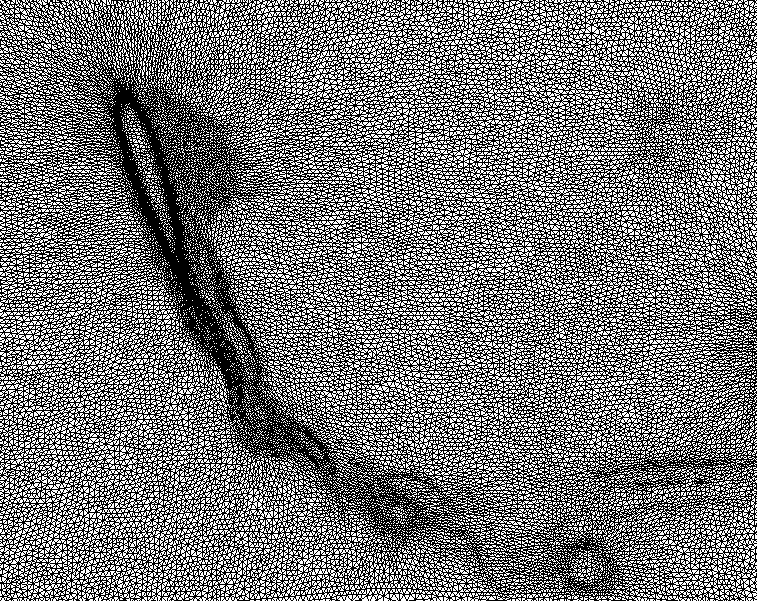}
\vspace{0.25cm}

\caption{Oscillation Naca snapshots. \textbf{top} Vorticity field. \textbf{bottom} Adaptive mesh.}
\label{fig:Naca_Sol}
\end{figure}

\begin{figure}
\centering
\includegraphics[width = 0.5\linewidth]{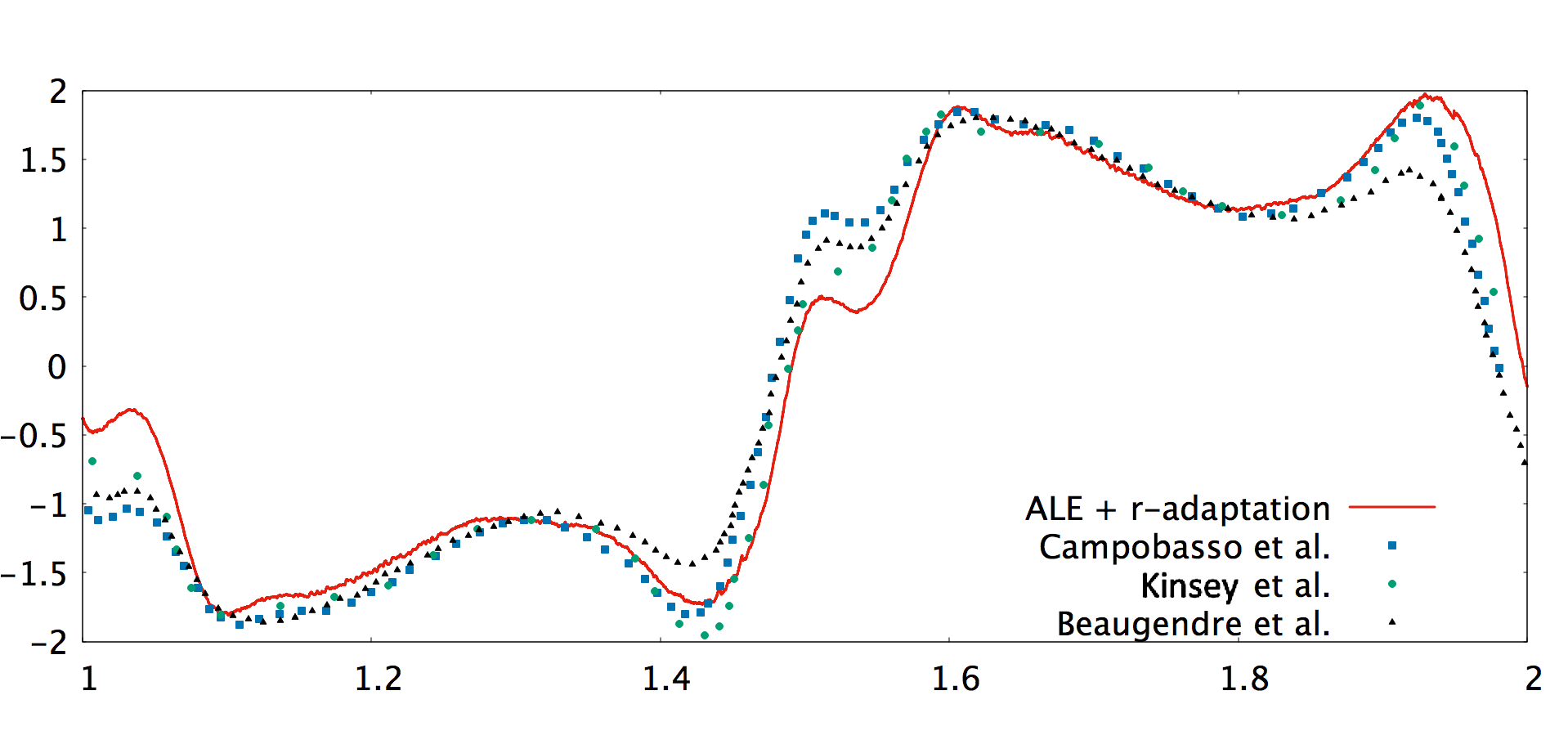}\includegraphics[width = 0.5\linewidth]{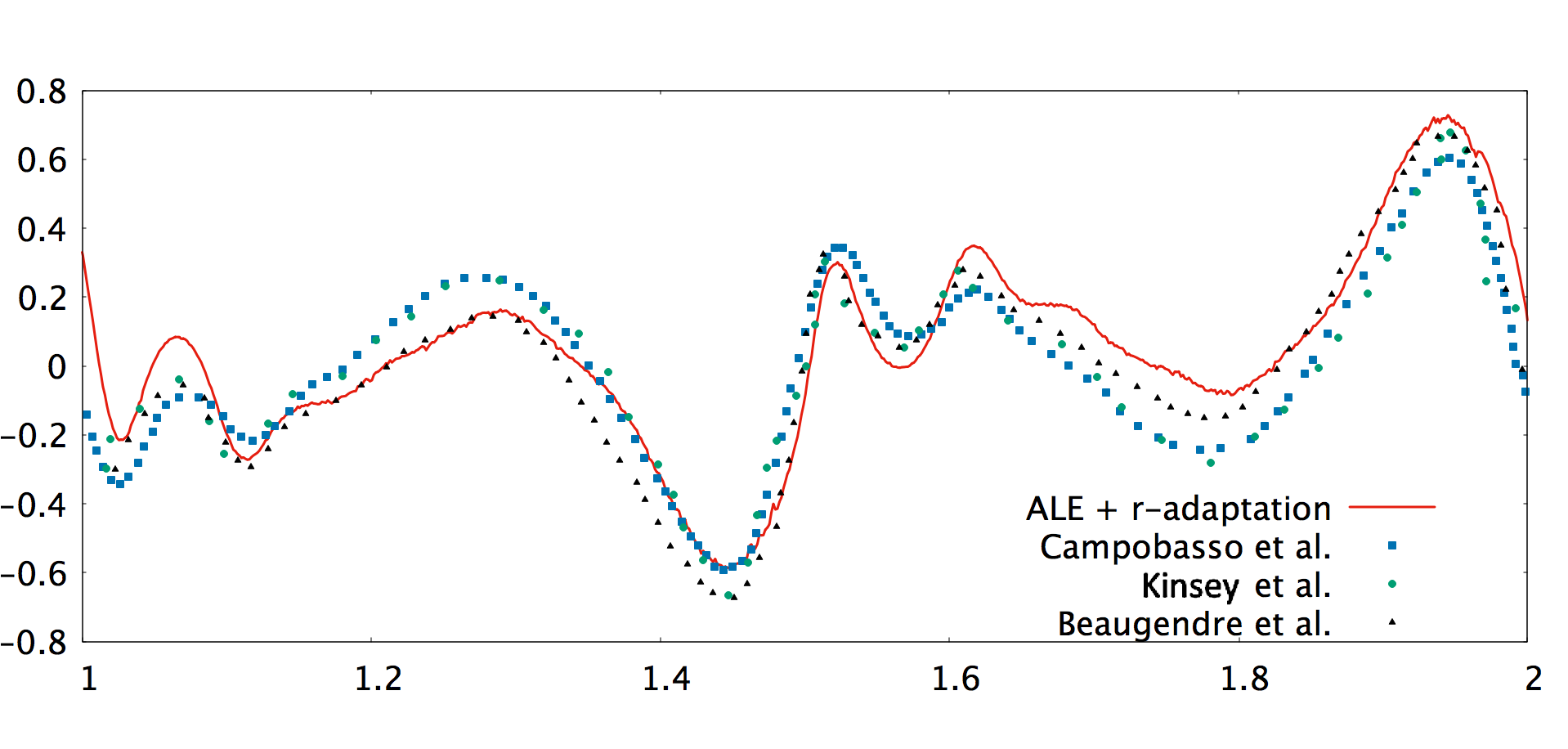}
\vspace{0.25cm}

\includegraphics[width = 0.5\linewidth]{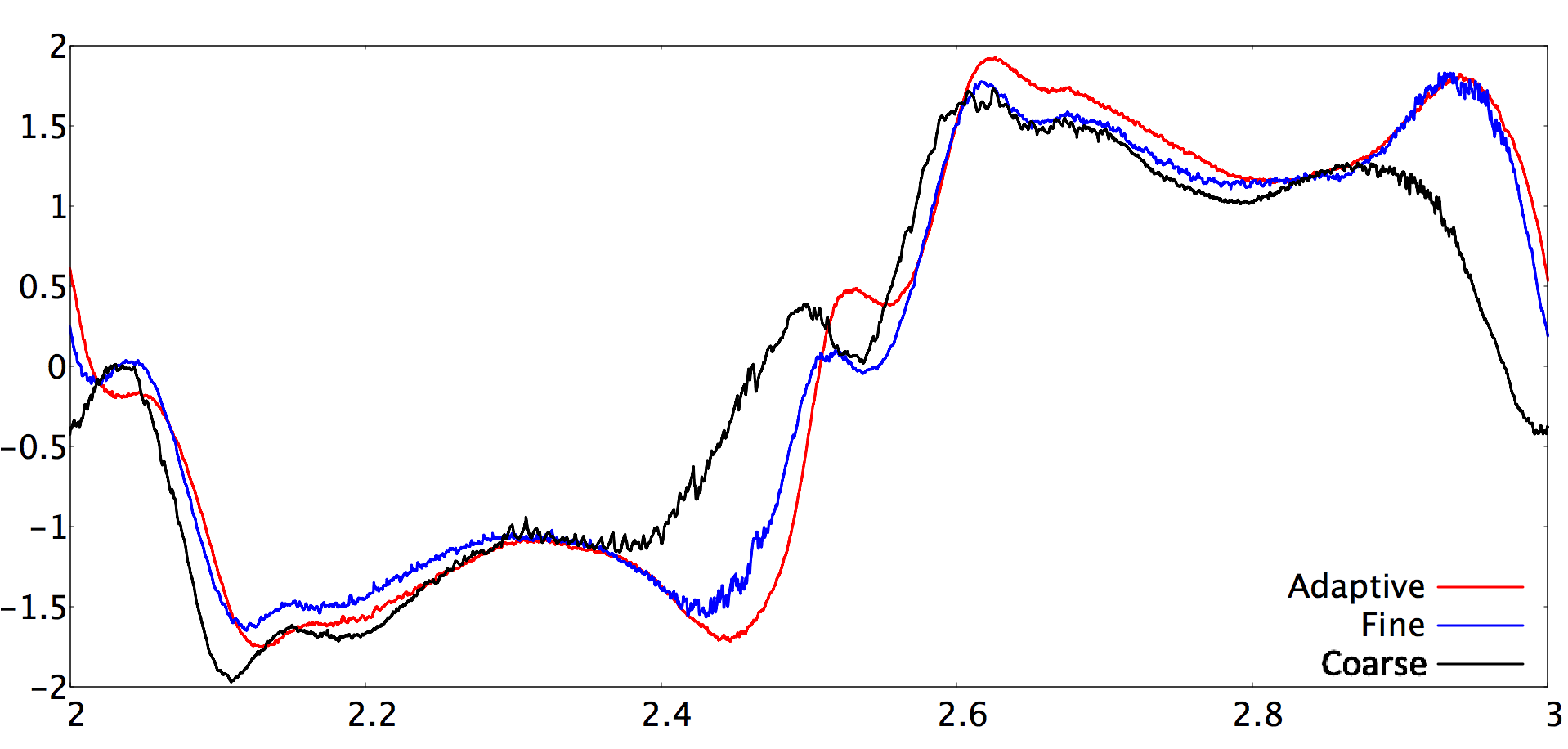}\includegraphics[width = 0.5\linewidth]{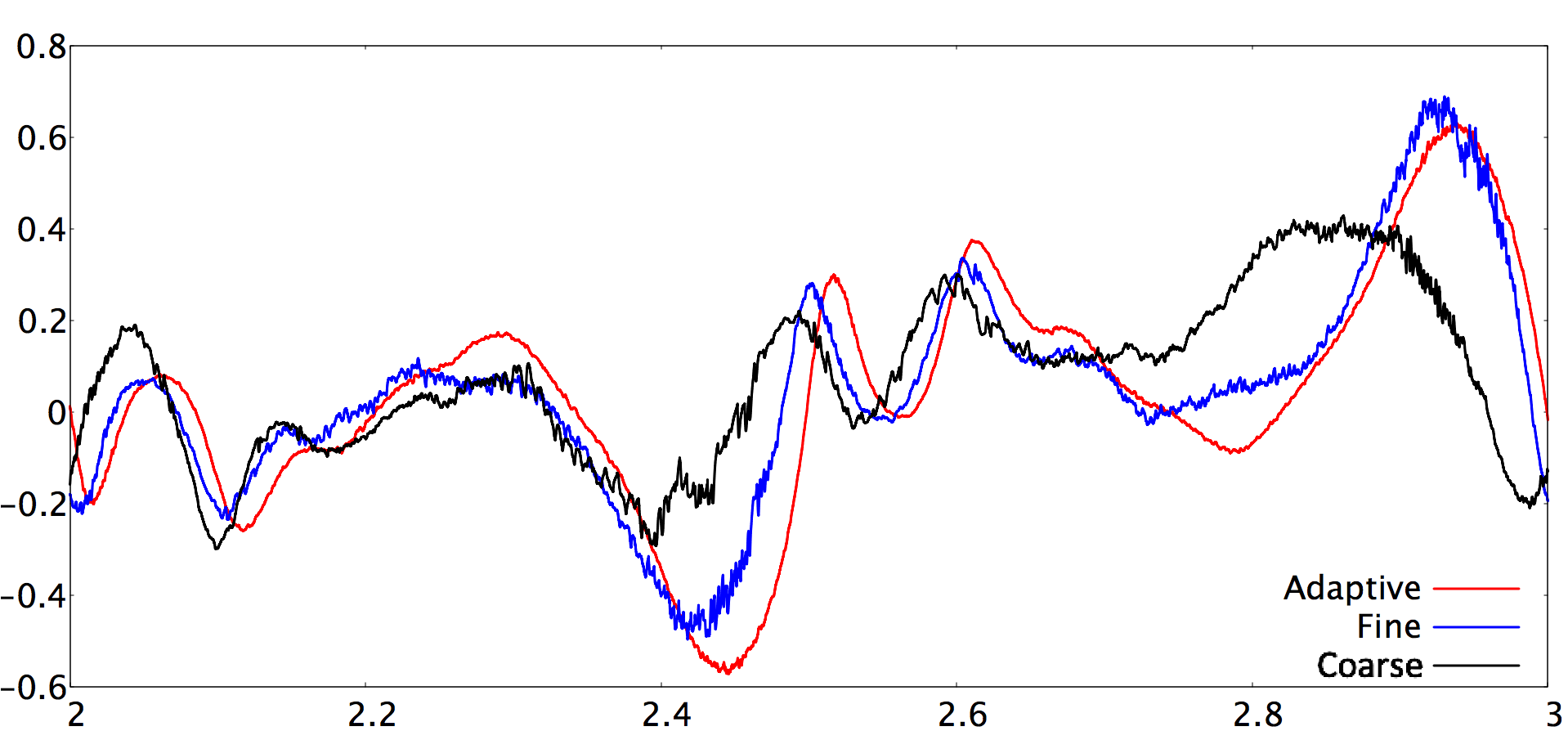}
\caption{Oscillating Naca - Aerodynamic coefficients vs time (Left: lift - Right: torque). \textbf{top} Comparison with literature. \textbf{bottom} Fine/adapted/coarse mesh results.}
\label{fig:Naca_Coeff}
\end{figure}

As an example we consider   a well known benchmark     involving the forced  pitching  and heaving  motion of a NACA 0015 airfoil
at  Reynolds $\approx1100$. We refer to \cite{Kinsey2008,Campobasso2012,Beaugendre2015}  for  the    set up  of the computation. 
Snapshots of the vorticity are reported in Fig.~\ref{fig:Naca_Sol}, showing the 
strong separation and vortex shedding taking place. In the same figure we report  the corresponding  meshes  adapted w.r.t. both the level-set, and vorticity. 
The efficient capturing of the body movement and vorticity fields in the mesh deformation is quite impressive. To validate our approach we compare the 
lift and torque aerodynamic coefficients to reference literature data ~\cite{Kinsey2008,Campobasso2012,Beaugendre2015} in the top pictures of Fig.~\ref{fig:Naca_Coeff}, observing an  excellent agreement.
We also compare the results obtained on the adaptive mesh, with those obtained on the initial reference triangulation (finest mesh size $h\approx 0.03$), and on a   refined one  (finest mesh size $h\approx 0.02$). 
Mesh data, and computational times are reported in Table~\ref{table:NacaTime}.  We can see that the $r-$adaptive approach allows computational savings of 
the order of   $26.5\%$,  and that without any particular optimization, only  $11.7\%$ of the computing time is devoted to the mesh PDE solver.

\begin{table}[!t]
\begin{center}
\caption{Oscillating Naca: computational times for different simulations}
\label{table:NacaTime}
\begin{tabular}{p{2.5cm}p{2.5cm}p{3cm}p{3.cm}}
\hline\noalign{\smallskip}
Mesh (size) & Number of nodes   & MMPDE CPU time [s]   & Total CPU time [s]\\
\noalign{\smallskip}\hline\noalign{\smallskip}
Fixed $h_{\min}\approx0.03$  & 53495 & $0$ & $20515$\\
Fixed $h_{\min}\approx0.02$ & 33139  &$0$ & $37716$\\
Adaptive &33139 & $3241$ & $27738$\\
\noalign{\smallskip}\hline\noalign{\smallskip}
\end{tabular}
\end{center}
\end{table}

\subsection{Free surface flows   and moving shorelines}

Another application with very similar issues is the simulation of free surface waves with moving wet/dry fronts.
These can be modelled by the shallow water equations which can be written as a system of balance equations
of the form
\begin{equation}
  \frac{\partial \mathbf{u}}{\partial t}
  +\nabla \cdot \mathbf{F}_C(\mathbf{u})
  + S( \mathbf{u}; \xbf) = {\bf 0}\;\;,\label{SWcompact}
\end{equation}
where $\mathbf{u}$ is  the array of the water depth $H$, and horizontal mass flux vector, $ \mathbf{F}_C$ the conservative fluxes, and $ S( \mathbf{u}; \xbf)$ a source term modelling several  
effects  such as those of bottom topography,   friction,    Coriolis  forces, and others depending on the application.  For  applications involving flooding and complex wave propagation 
ideas similar to those discussed above can be used to  adapt the mesh w.r.t  shorelines and wave fronts.  However, for wave propagation, and in general for hyperbolic problems, 
a majority of codes are based on an  explicit time stepping kernel. This requires limiting somewhat the 
overhead of the mesh PDE solver  for $r$-adaptation to be interesting. 
So we do not   use  the computation of a distance function here, but  we  proceed as follows:
\begin{itemize}
\item generate initial meshes with sizes increasing with the depth, so that a prescribed minimum size is attained at and close to the initial shoreline
\item define the regularized Heaviside function  
$$
\phi_H(\xbf)  = \left\{\begin{array}{lll}
1 \quad&\text{if }\;  H(\xbf) > \epsilon_H\\[5pt]
\phi(H) \quad&\text{if }\;  0< H(\xbf)  \le \epsilon_H\\[5pt]
0 \quad&\text{if }\;  H(\xbf) = 0 
\end{array}
\right.
$$
\item define a  monitor function combining the dry detection function above and  a capped gradient  (cf. Eq.~\ref{eq:nablac-cap}) as 
\begin{equation}\label{eq:omega4}
\omega = \sqrt{1+ \alpha_{\eta}\widehat{\nabla \eta}^2 + \alpha_\text{dry}\|\nabla_{\xbf} \phi_H \|^2}
\end{equation} 
\end{itemize}
where if $b({\xbf})$ denotes the topography,  $\eta$ is the free surface level $H+b$. 
The  regularization $\phi(H)\in[0,1]$ can be any (at least $C^0$) function modelling the depth profile in cells cut by the shoreline. 
Note that, differently than in the case of moving bodies,
this curve is in general never known exactly.  Even in the initial condition, its definition depends in practice on the resolution of the topographic data.
In practice,  defining $\phi(H)$ as a  linear function is already enough to capture the wet/dry transition.   \\

In the numerical examples shown here, the shallow water equations are solved by means of an  explicit stabilized
residual based  approach  in  an ad-hoc ALE form. Although the design of numerical schemes is not the focus of this  work,
we mention that one of the challenges  when dealing with the shallow water equations (and in general balance laws)
 is that the ALE formulation, and the associated remaps, needs to satisfy not only the classical conservation  constraints 
(mass and momentum), but also additional ones related to the preservation of particular steady solutions (well-balanced).
As for compressible fluid flows, additional constraints related to the admissibility of the solution (e.g.  non-negativity of the depth) should also be embedded. 
Moreover, in real applications one has to make sure that  the remap used for the data of the problem  (e.g. topography) is not 
 affected by truncation errors, which is also a constraining requirement.
These aspects are covered in detail in \cite{arpaia:hal-01372496}, in which these constraints are 
 combined within a conservative ALE formulation, with a high order  quadrature based projection of  the original data on the moving mesh.
 
 Concerning the mesh solver, while for the initial solution  Eq.~\ref{eq:Jacobi0} is solved until convergence within a few orders of  magnitude,
 within each time step only 5 Jacobi relaxation iterations are performed to move the mesh, using as initial guess the last available one.
As an example, we consider the simulation of the  2011 Tohoku-Honsu tsunami, the interested reader can refer to \cite{arpaia:hal-01372496,Arpaia2020}
for thorough validation.

\begin{figure}
\includegraphics[width=0.27\textwidth]{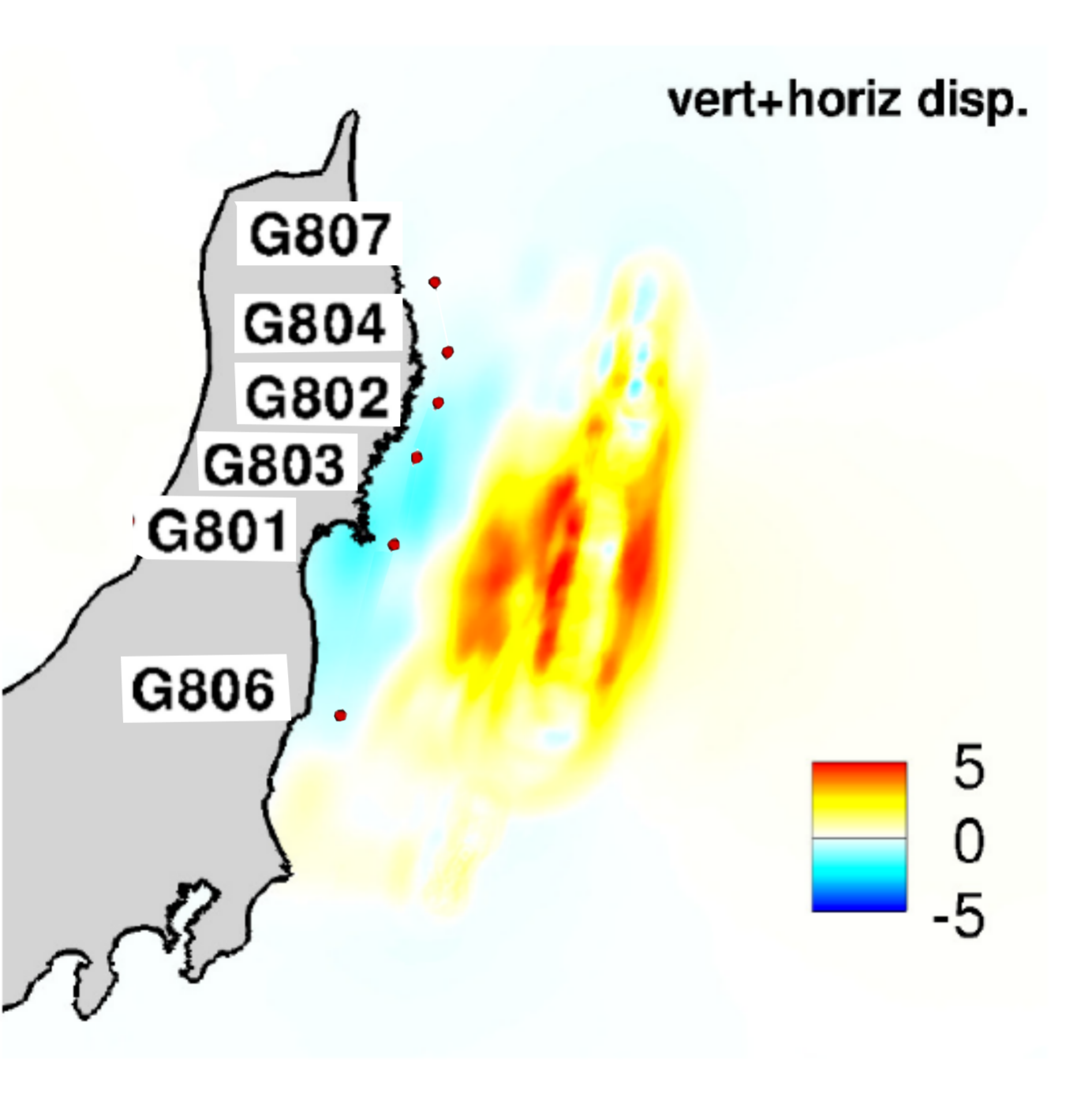}\includegraphics[width=0.4\textwidth]{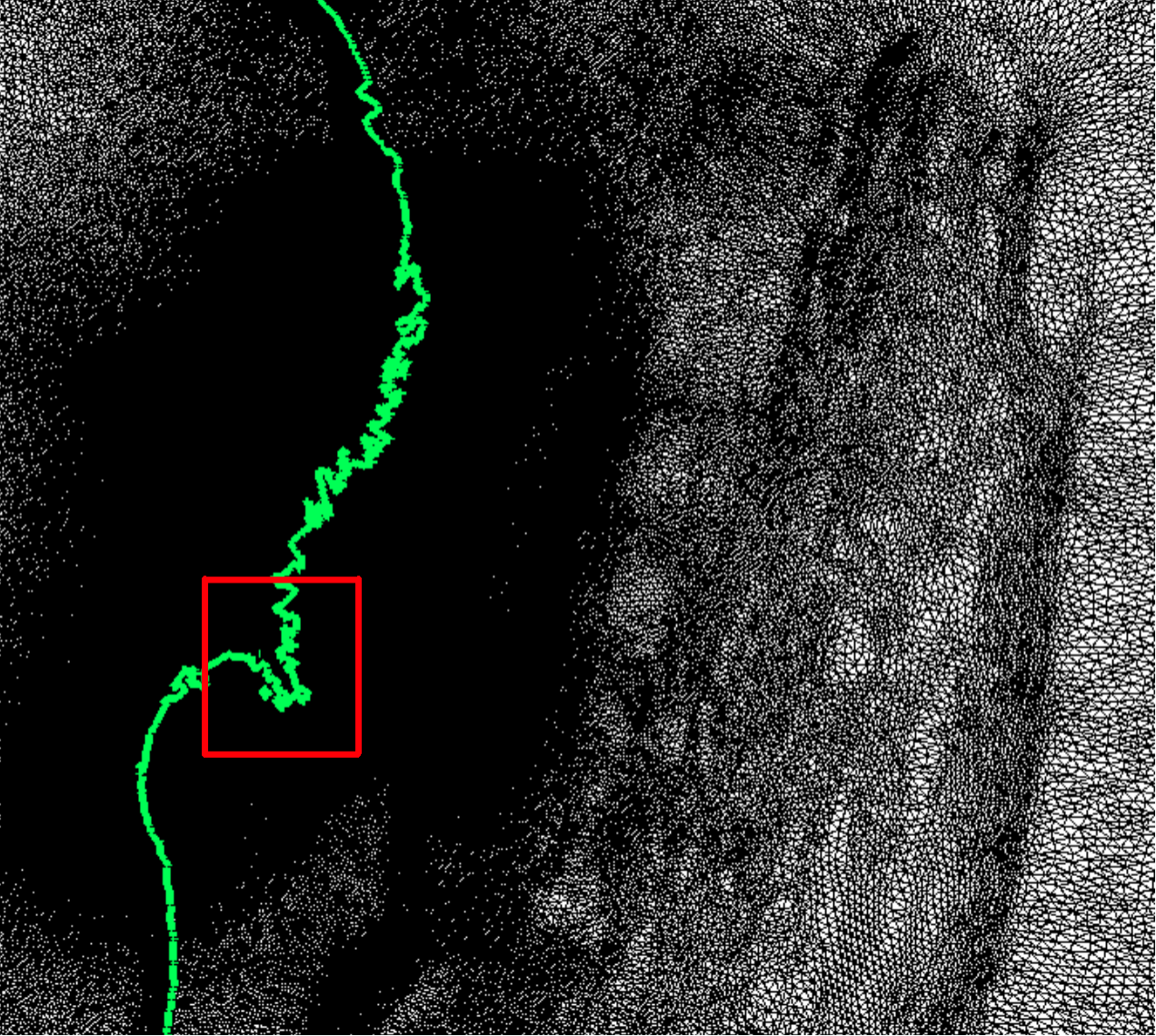}\includegraphics[width=.3\textwidth]{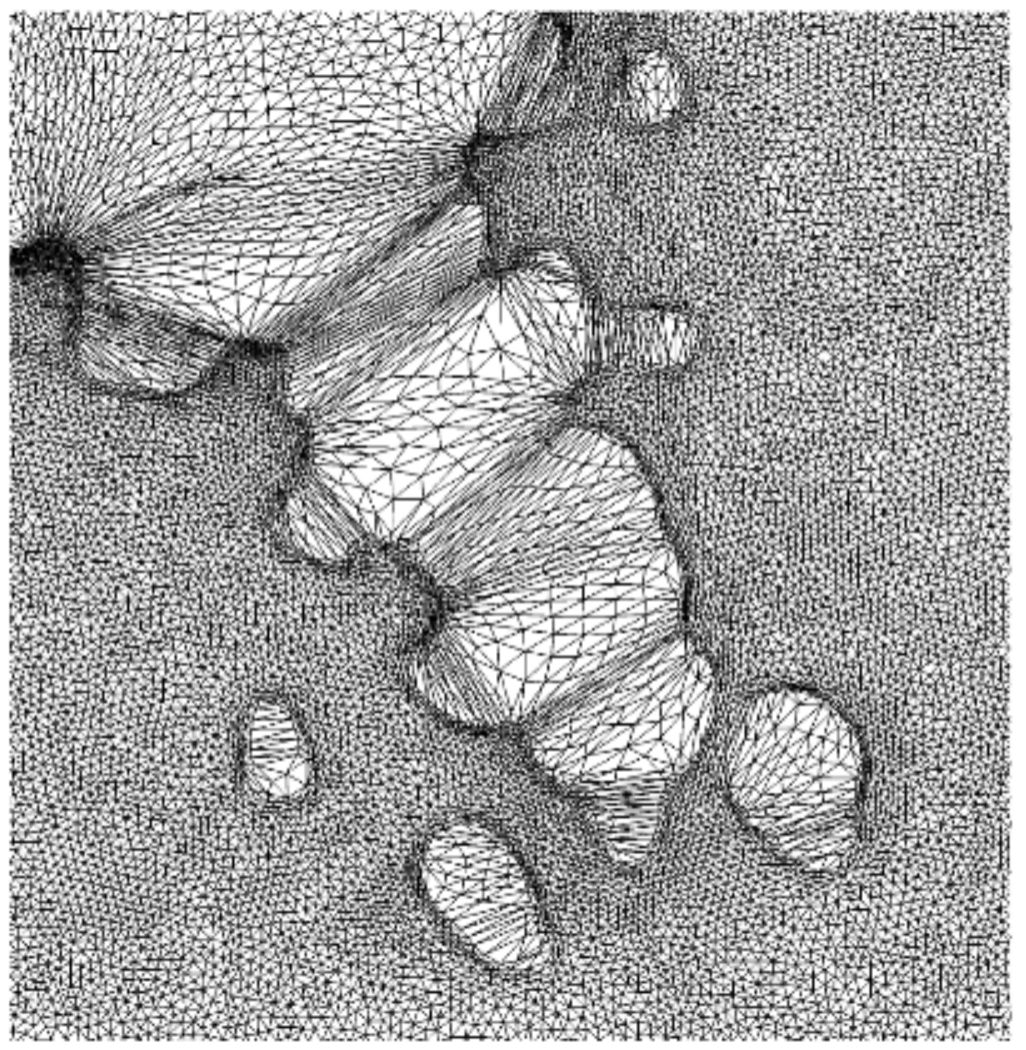}
\caption{ Tohoku-Honsu tsunami. \textbf{left} Problem setup, GPS probes and initial wave height (in meters). \textbf{center} Adapted closer to the coastline
of the Sendai bay and Myagi prefecture (boxed area). \textbf{right} Shoreline of the Myagi prefecture.}
\label{fig:sendai-1}
\end{figure}

The leftmost picture on Fig.~\ref{fig:sendai-1} shows an overview of the computational domain, including the position of the GPS buoys whose data has been used for validation,
and  the  initial free surface perturbation   computed  based on the seismic data and  with approach by Satake \cite{satake,brgm17}.
The middle picture in the same figure shows a close up of the Japanese coast in correspondence of the Myagi prefecture (boxed area) and of the Sendai bay, with
the adapted mesh in the initial instants of the  propagation of the tsunami.  
The rightmost picture shows a further close up of the Myagi prefecture  with the mesh adapted to the coastline. Note that the resolution of the  reference mesh goes from 360 m close to the initial shoreline to 15 Km offshore. 
A fine mesh computation has been run for comparison using resolutions  down to 120m close to the shore and of 5 Km offshore.

\begin{figure}
\begin{minipage}{0.62\textwidth}
\includegraphics[width=0.5\textwidth]{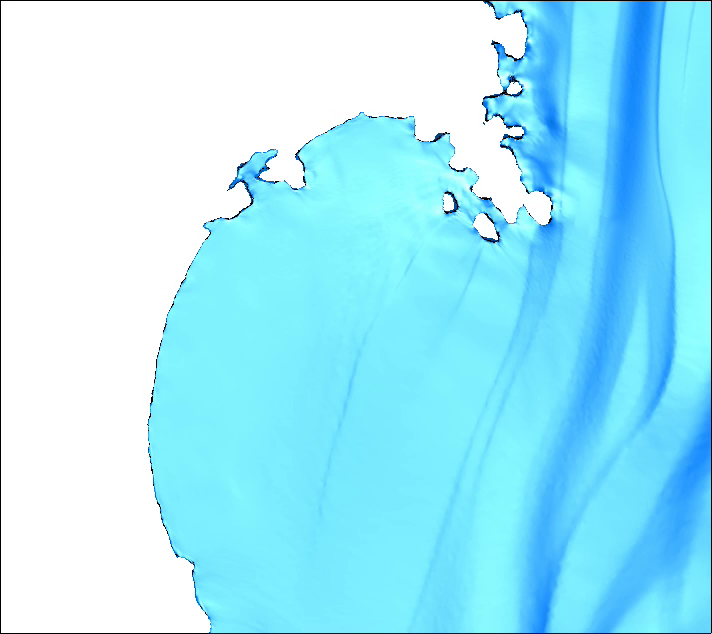}
\includegraphics[width=0.5\textwidth]{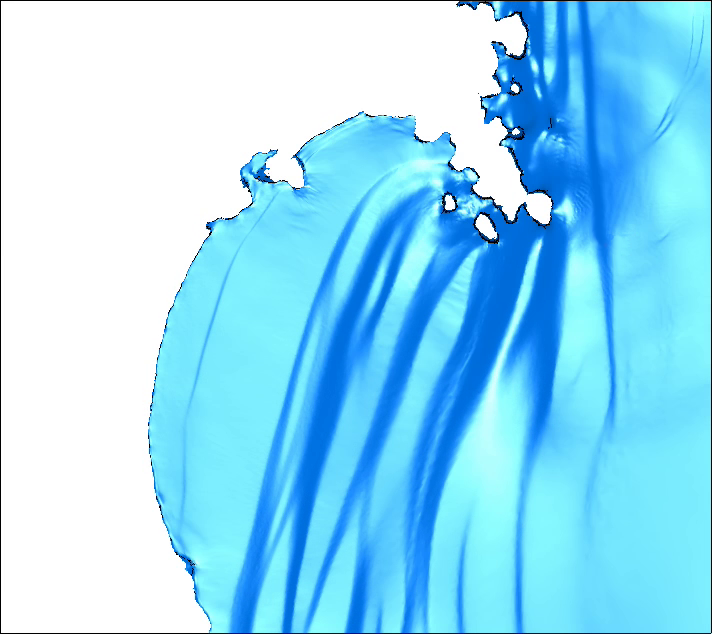}\\[5pt]

\includegraphics[width=0.5\textwidth]{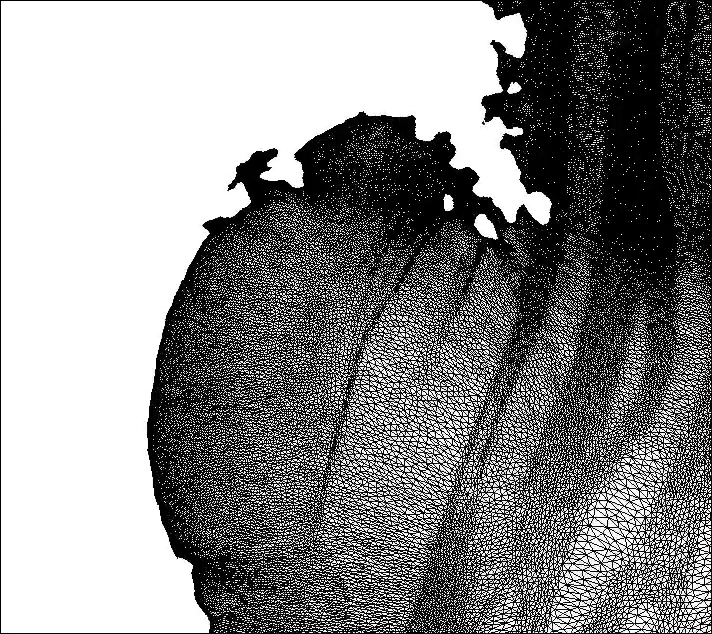}
\includegraphics[width=0.5\textwidth]{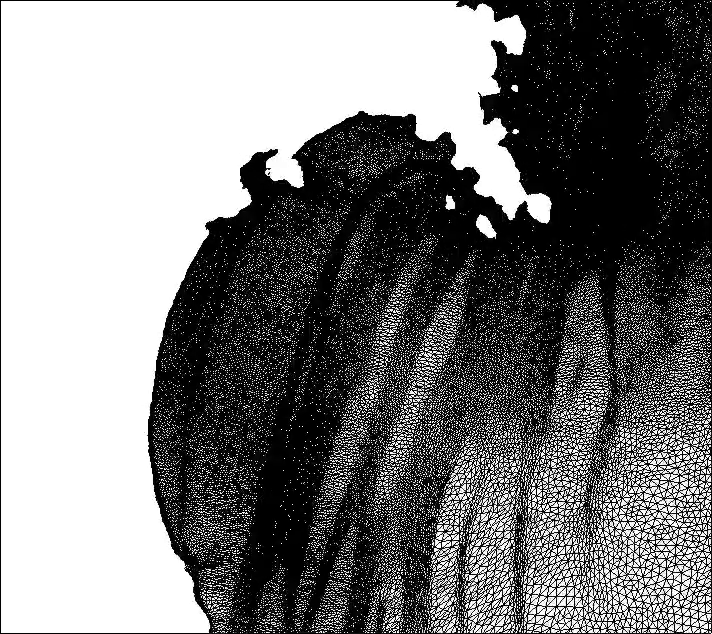}
\end{minipage}
\begin{minipage}{0.375\textwidth}
\includegraphics[width=1.\textwidth]{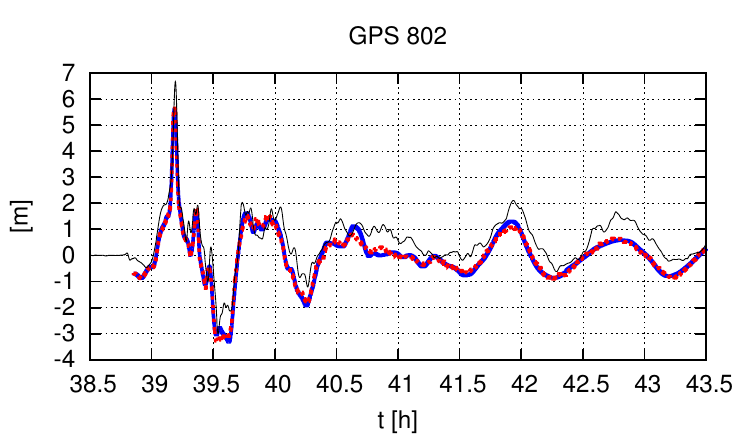}\\[10pt]
\includegraphics[width= 1.\textwidth]{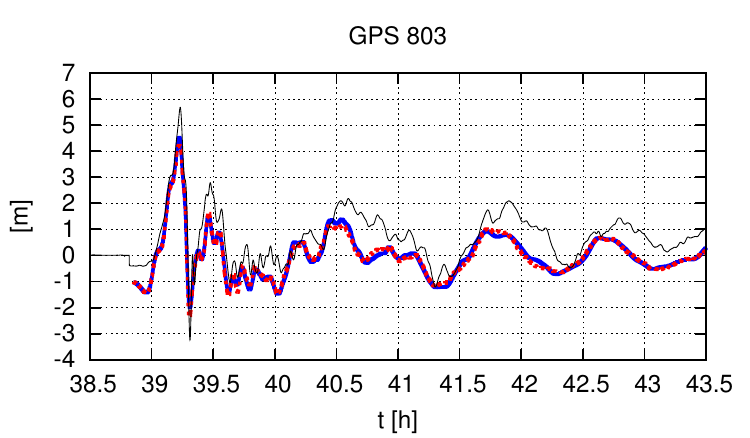}
\end{minipage}
\caption{ Tohoku-Honsu tsunami. \textbf{left and center} Waves entering the Sendai bay. \textbf{right} GPS bouys 802 and 803 (red: adaptive - blue:  fine mesh - black: GPS data).}
\label{fig:sendai0}
\end{figure}

\begin{figure}
\centering
\includegraphics[width=0.35\textwidth,trim=1mm 1mm 1mm 1mm,clip]{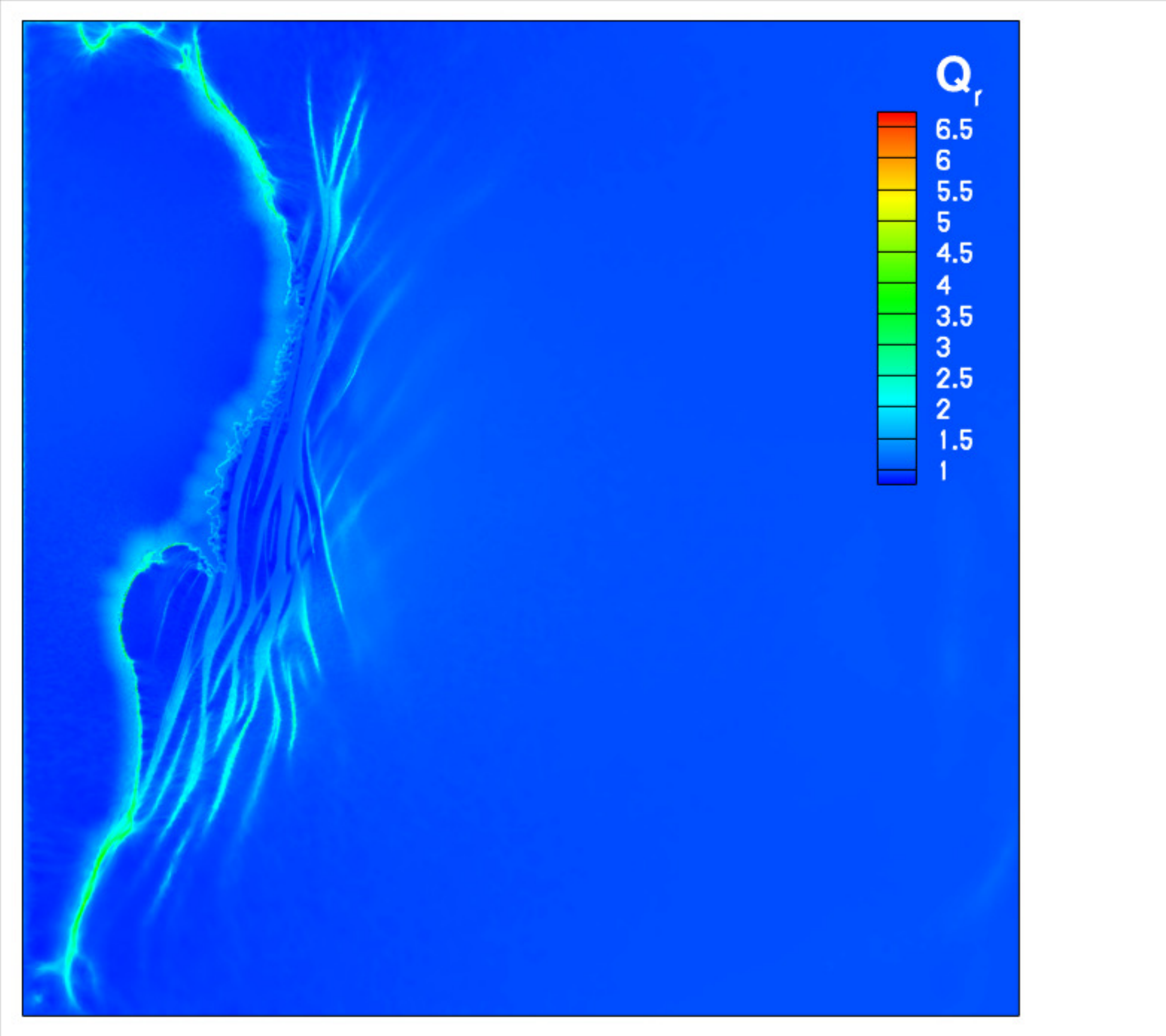}\hspace{0.4cm}
\includegraphics[width=0.35\textwidth,trim=1mm 1mm 1mm 1mm,clip]{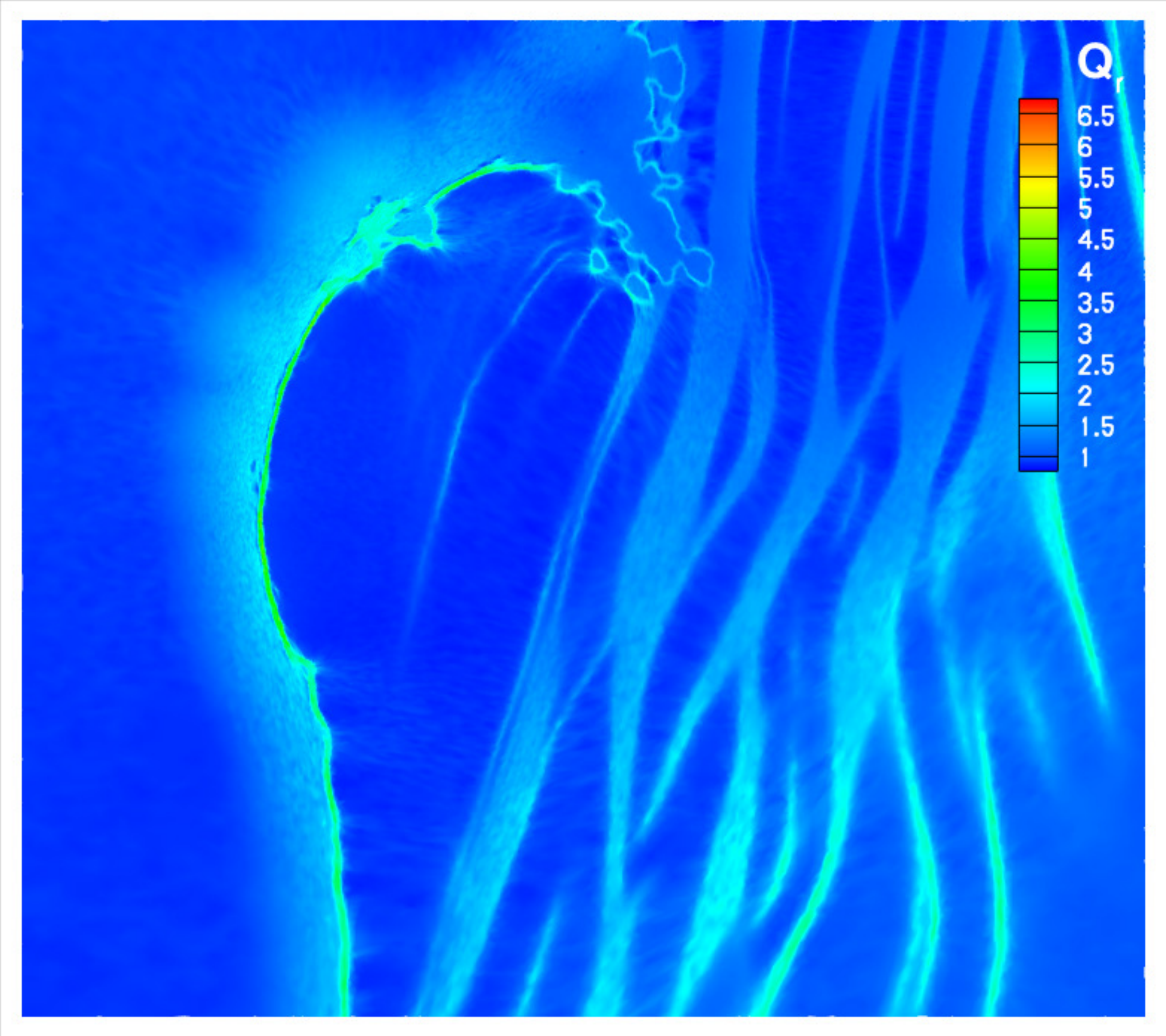}\\[3pt]
\includegraphics[width=0.35\textwidth,trim=1mm 1mm 1mm 1mm,clip]{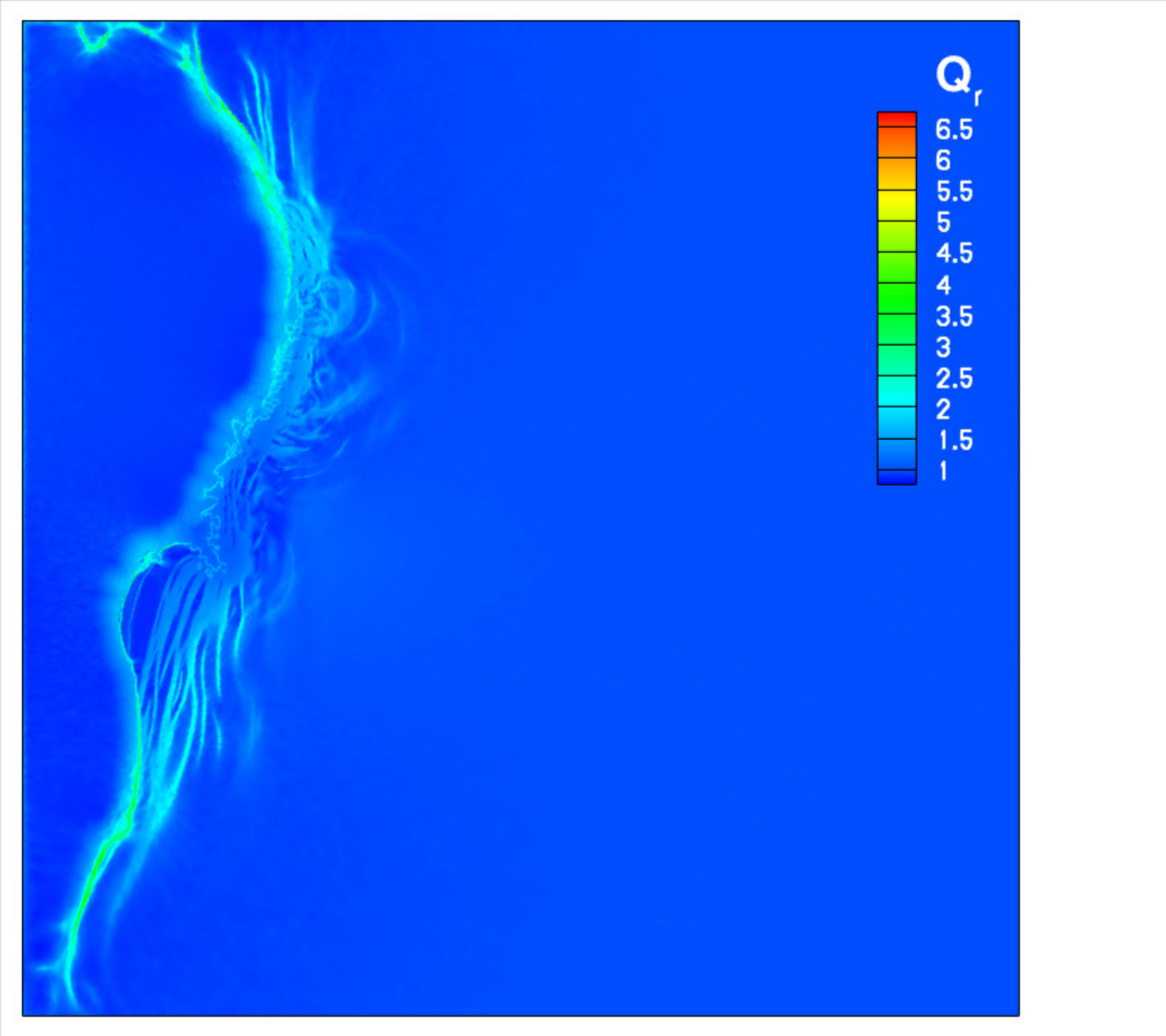}\hspace{0.4cm}
\includegraphics[width=0.35\textwidth,trim=1mm 1mm 1mm 1mm,clip]{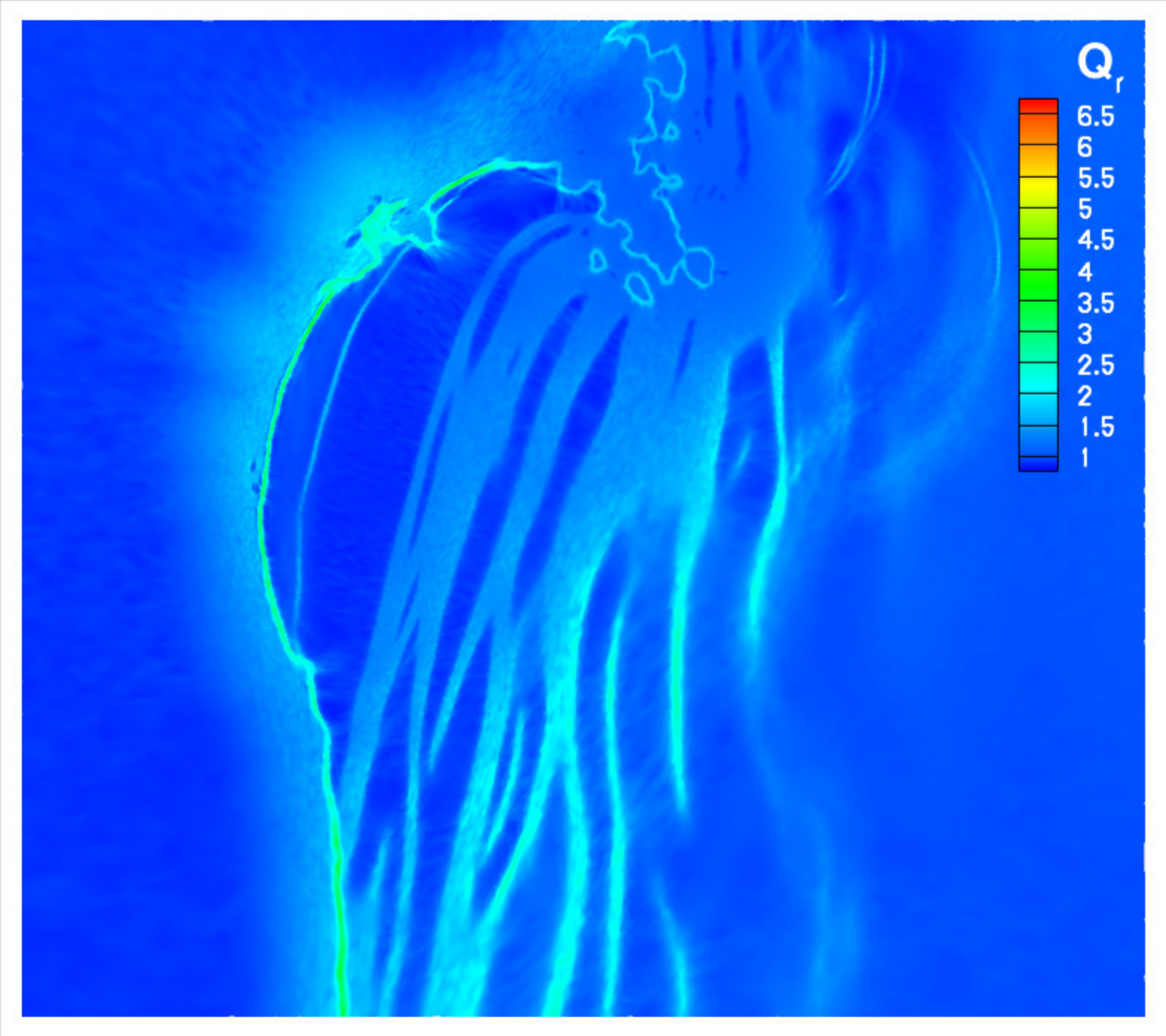}
    \caption{Tohoku-Honsu tsunami. Anisotropic compression, $\Qcal_r$, of the adapted mesh. \textbf{left} Whole computational domain. \textbf{right} Close up of the Sendai bay.}\label{fig:compression_tsunami}  
\end{figure}

The left and middle pictures on Fig.~\ref{fig:sendai0} show the waves entering the Sendai bay, and the corresponding adapted meshes,
while the rightmost pictures report the GPS time series, comparing the adaptive computations to the fine mesh, and to the GPS data showing the potential
of the approach used in resolving complex wave propagation on coarser grids. \mario{Finally, Fig.~\ref{fig:compression_tsunami} shows the distribution of    $\Qcal_r$, defined in Eq.~\ref{eq:compdef}, at two different time steps with a close up of the Sendai bay.  The compression ratios achieved are of about  $2.5-3$  across the steepest shoaling waves, and increase above $4$ and with peaks  of  about $6-6.5$ in correspondence of the shorelines.}

